\documentclass[11pt]{article}
\usepackage{eurosym}
\usepackage{amsfonts}
\usepackage{amsmath}
\usepackage{geometry}
\usepackage{amssymb}
\usepackage{array}
\usepackage{graphicx}
\usepackage{color}
\usepackage[dvipsnames]{xcolor}
\usepackage{mathrsfs}
\usepackage[author={Oana}]{pdfcomment}
\usepackage[colorinlistoftodos]{todonotes}
\usepackage{titling}
\usepackage{lipsum}
\usepackage{soul}
\usepackage{dsfont}

\newtheorem{theorem}{Theorem}

\newtheorem{corollary}[theorem]{Corollary}

\newtheorem{definition}[theorem]{Definition}

\newtheorem{lemma}[theorem]{Lemma}

\newtheorem{proposition}[theorem]{Proposition}
\newtheorem{remark}[theorem]{Remark}

\newcommand{\T}{\mathbb{T}}

\newcommand{\cL}{\mathcal{L}}
\newcommand{\cA}{\mathcal{A}}
\newcommand{\cT}{\mathcal{T}}
\newcommand{\cW}{\mathcal{W}}

\newcommand{\vertiii}[1]{{\left\vert\kern-0.25ex\left\vert\kern-0.25ex\left\vert #1 
    \right\vert\kern-0.25ex\right\vert\kern-0.25ex\right\vert}}

\begin{document}

\date{}
\title{\textbf{Well-posedness Properties for a Stochastic Rotating Shallow
Water Model}}
\author{Dan Crisan \qquad Oana Lang \\
\\
{\small Department of Mathematics, Imperial College London, SW7 2AZ, UK}}
\maketitle

\begin{abstract}
In this paper, we study the well-posedness properties of a stochastic rotating
shallow water system. An inviscid version of this model has first been derived in \cite{Holm2015} and the noise is chosen according to the Stochastic Advection by Lie Transport theory presented in \cite{Holm2015}. 
The system is perturbed by noise modulated by a
function that is not Lipschitz in the norm where the
well-posedness is sought. We show
that the system admits a unique maximal solution which depends
continuously on the initial condition. We also show that the interval of existence is strictly positive and the solution is global with positive probability. 

\end{abstract}

\tableofcontents

\section{Introduction}

The rotating shallow water equations describe the evolution of a
compressible rotating fluid below a free surface. The typical vertical
length scale is assumed to be much smaller than the horizontal one, hence
the \textit{shallow} aspect. This model is a simplification of the primitive
equations which are known for their complicated and computationally
expensive structure (see e.g. \cite{Kalnay}). Despite its simplified form, the
rotating shallow water system retains key aspects of the atmospheric and
oceanic dynamics (\cite{Kalnay}, \cite{Vallis}, \cite{Zeitlinbook}). It allows
for gravity waves which play a highly important role in climate and weather
modelling (\cite{Vallis}). The classical inviscid shallow water model
consists of a horizontal momentum equation and a mass continuity equation
and in the presence of rotation it can be described as follows (see \cite%
{Darryllecturenotes}):

\begin{equation*}
\epsilon \frac{D}{Dt}u_t + f\hat{z}\times u_{t}+\nabla p_{t} =0 
\end{equation*}
\begin{equation*}
\frac{\partial h_{t}}{\partial t}+\nabla \cdot (h_{t}u_{t}) =0
\end{equation*}
where
\begin{itemize}
\item $\frac{D}{Dt} := \frac{\partial}{\partial t} + u \cdot \nabla$ is the
material derivative

\item $u=(u^1,u^2)$ is the horizontal fluid velocity vector field 
\item $\epsilon$ is the Rossby number, a dimensionless number which
describes the effects of rotation on the fluid flow: a small Rossby number ($%
\epsilon <<1$) suggests that the rotation dominates over the advective
terms; it can be expressed as $\epsilon = \frac{U}{fL}$ where $U$ is a
typical scale for horizontal speed and $L$ is a typical length scale.

\item $f$ is the Coriolis parameter, $f=2\Theta \sin\varphi$ where $\Theta$
is the rotation rate of the Earth and $\varphi$ is the latitude; $f\hat{z}
\times u = (-fu^2, fu^1)$, where $\hat{z}$ is a unit vector pointing away
form the centre of the Earth. For the analytical analysis we assume $f$ to
be constant.
\item $h$ is the thickness of the fluid column
\item $p := \frac{h-b}{\epsilon\mathscr{F}}$, $\nabla p$ is the pressure
term, $b$ is the bottom topography function.
\item $\mathscr{F}$ is the Froude, a dimensionless number which relates to
the stratification of the flow. It can be expressed as $\mathscr{F}=\frac{U}{%
NH}$ where $H$ is the typical vertical scale and $N$ is the buoyancy
frequency.
\end{itemize}

The deterministic nonlinear shallow water equations (also known as the 
\textit{Saint-Venant equations}) have been extensively studied in the
literature. A significant difficulty in the well-posedness analysis of this
model is generated by the interplay between its intrinsic nonlinearities, in
the absence of any incompressibility conditions. In order to counterbalance
the resulting chaotic effects, a viscous higher-order term is usually added
to the inviscid system. Various shallow water models have been introduced for instance in
\cite{Zeitlinbook} or \cite{BreschsDesjardinsMetivier}. In \cite{LiuYin2} the
authors show global existence and local well-posedness for the 2D viscous
shallow-water system in the Sobolev space $H^{s-\alpha}(\mathbb{R}^2) \times
H^s(\mathbb{R}^2)$ with $s>1$ and $\alpha \in [0,1)$. The methodology is
based on Littlewood-Paley approximations and Bony paraproduct
decompositions. This extends the result in \cite{LiuYin1} where local solutions
for any initial data and global solutions for small initial data have been
obtained in $H^s \times H^s$ with $s>1$. A similar result adapted to Besov
spaces was obtained in \cite{LiuYin0}. More recently, ill-posedness for the
two-dimensional shallow water equations in critical Besov spaces has been
shown in \cite{LiHongZhu}. Existence of global weak solutions and convergence to
the strong solution of the viscous quasi-geostrophic equation, on the
two-dimensional torus is shown in \cite{BreschsDesjardins1}. In
\cite{BreschsDesjardins2} the authors construct a sequence of smooth approximate
solutions for the shalow water model obtained in \cite{BreschsDesjardins1}. The
approximated system is proven to be globally well-posed, with height bounded
away from zero. Global existence of weak solutions is then obtained using
the stability arguments from \cite{BreschsDesjardins1}. Sundbye in \cite{Sundbye}
obtains global existence and uniqueness of strong solutions for the
initial-boundary-value problem with Dirichlet boundary conditions and small
forcing and initial data. In this work the solution is shown to be classical
for a strictly positive time and a $C^0$ decay rate is provided. The proof
is based on a priori energy estimates. Independently, Kloeden has shown in
\cite{Kloeden} that the Dirichlet problem admits a global unique and spatially
periodic classical solution. Both \cite{Sundbye} and \cite{Kloeden} are based on the
energy method developed by Matsumura and Nishida in \cite{Matsumura}. Local
existence and uniqueness of classical solutions for the Dirichlet problem
associated with the non-rotating viscous shallow water model with initial
conditions $(u_0,h_0) \in C^{2, \alpha} \times C^{1, \alpha}$ can be found
in \cite{Bui}. The proof is based on the method of successive approximations and Hölder space estimates, in a Lagrangian framework. Existence and
uniqueness of solutions for the two-dimensional viscous shallow water system
under minimal regularity assumptions for the initial data and with height
bounded away from zero was proven in \cite{ChenMiaoZhang}. The possibly
stabilising effects of the rotation in the inviscid case is analysed in
\cite{rotRSW1} and \cite{rotRSW2}.

To simplify notation, we will denote by $a:=\left( v,h\right) $ the solution
of the rotating shallow water (RSW) system and recast it in short form as\footnote{%
We use here the differential notation to match the stochastic version (\ref%
{shortsrsw}).}
\begin{equation*}
da_{t}+F\left( a_{t}\right)dt =0,
\end{equation*}%
where $F\left( a_{t}\right) $ denotes 
\begin{equation*}
F\left( 
\begin{array}{c}
v \\ 
h%
\end{array}
\right) =\left( 
\begin{array}{c}
u \cdot \nabla v + f\hat{z} \times u + \nabla p \\ 
\nabla \cdot (hu)
\end{array}
\right)
\end{equation*}
where $u$ is the fluid velocity and $v: =\epsilon u + \mathcal{R}$, with $curl \ \mathcal{R}
= f\hat{z}$. $\mathcal{R}$ corresponds to the vector potential for the
(divergence-free) rotation rate about the vertical direction, and it is
chosen here such that $\nabla \mathcal{R} = 0$.
In this paper we consider a viscous and stochastic version of the shallow
water model described above, defined on the two-dimensional torus $\mathbb{T}%
^{2}$:

\begin{equation}  \label{shortsrsw}
da_{t}+F\left( a_{t}\right)dt +\sum_{i=1}^{\infty }\mathcal{G}_{i}\left(
a_{t}\right) \circ dW_{t}^{i}=\gamma \Delta a_{t}dt,
\end{equation}%
where $\gamma = (\nu,\eta) $ is positive and corresponds to the fluid viscosity\footnote{%
Different levels of viscosity for the different components of $a$ can be
treated in the same manner.}, $W^{i}$ are independent Brownian motions, $F$
is a nonlinear advective term, and $\mathcal{G}_{i}$ are differential
operators explicitly described below. The integrals in (\ref{shortsrsw}) are
of Stratonovitch type. The system (\ref{shortsrsw}) belongs to a class of
stochastic models derived using the Stochastic Advection by Lie Transport
Approach (SALT) approach, as described in \cite{Holm2015}, \cite{SC2020}, \cite%
{HolmLuesink}. 
A detailed derivation of
this specific system can be found in \cite{OanaPhD}, following \cite%
{HolmLuesink}, \cite{Holm2015}. In the stochastic case, $F$ is defined as
above and 
\begin{equation*}
\mathcal{G}_{i}\left( 
\begin{array}{c}
v \\ 
h%
\end{array}%
\right) =\left( 
\begin{array}{c}
\cL_i v + \cA_i v \\ 
\cL_i h%
\end{array}%
\right) 
\end{equation*}%
where $\xi _{i}$ are divergence-free and time-independent vector fields,
$\mathcal{L}_iv = \xi_i \cdot \nabla v$, $\mathcal{A}_iv = v^1\nabla\xi_i^1+v^2\nabla\xi_i^2$. The two operators $\mathcal{L}_i$ and $\mathcal{A}_i$ enjoy some
properties which are described in Section \ref{preliminariesSRSW} and in the Appendix. 
It
has been shown lately that such stochastic parameters can be calibrated
using data-driven approaches to account for the missing small-scale
uncertainties which are usually present in the classical deterministic
geophysical fluid dynamics models (see for instance \cite{CotterCrisanPan1},
\cite{CotterCrisanPan2}). The addition of stochasticity in the advective part
of the dynamics brings forth a more explicit representation of the uncertain
transport behaviour in fluids, which draws on recent synergic advances in
stochastic analysis, geophysical fluid dynamics, and data assimilation. The
performance of these modern stochastic approaches is subject of
intensive research. In a forthcoming work
(\cite{LeeuwenCrisanPotthastLang}) we prove the applicability of this new
stochastic model in a data assimilation framework.

To the best of our knowledge, this specific form of the stochastic rotating
shallow water model has not been studied before. A stochastic version of the
viscous rotating shallow water system with external forcing and
multiplicative noise has been studied in \cite{Temam1}. This corresponds to
the case $\mathcal{G}_{i}\left( a\right) =a$. 
A rotating shallow water model driven by L\'evy noise has been considered in 
\cite{Temam2}. As pointed out in \cite{Temam1}, the number of results
available in the literature on stochastic shallow water equations is
limited. In the deterministic case, existence of solutions under certain
conditions and without rotation was proven in \cite{Orenga}. Smooth
approximate solutions for the 2D deterministic rotating shallow water system
have been constructed in \cite{Breschs2}. Long time existence for rapidly
rotating deterministic shallow water and Euler equations has been shown in 
\cite{rotRSW1}.

\subsection{Contributions of the paper}

The first contribution of the paper is the existence of local solutions $%
(a,\tau)$ of the system (\ref{shortsrsw}) with paths $t\rightarrow
a_{t\wedge\tau}$ in the space\footnote{$\mathcal{W} ^{k,2}$ is the standard
Sobolev space.} 
\begin{equation*}
C\left( [0,T],\left( \mathcal{W}^{1,2}(\mathbb{T}^2)\right) ^{3}\right) \cap
L^{2}\left(0,T;\left( \mathcal{W}^{2,2}\right) ^{3}(\mathbb{T}^2)\right) ,~\
~T>0
\end{equation*}%
provided the initial datum $a_{0}\in \left( \mathcal{W}^{1,2}\right) ^{3}\ $%
($\tau $ is a strictly positive stopping time to be specified below).
Subsequently, we prove that there exists a unique strong maximal solution $%
(a,\tau_{\mathrm{max }})$ of the system (\ref{shortsrsw}), see Theorem \ref%
{mainthm} below. We also show that the solution depends continuously on the
initial data, see Theorem \ref{prop:uniquenesstruncated} below. 

The approach we follow uses a priori estimates similar to those in \cite%
{Temam1}, as the deterministic terms are very similar. However, the structure of the noise studied here is different
from that in \cite{Temam1}, in particular the operators ${\mathcal{G}}_i$
are not Lipschitz. As a result, we need to use a different approximation
method, extending the methodology developed in \cite{CL1} and \cite{CL2} to a system
of SPDEs which model a compressible fluid under the effects of rotation. In
particular, we do not use Galerkin approximations as in \cite{Temam1}.
Instead we construct a sequence that approximates a truncated version of the equation, where
the nonlinear term is replaced by a forcing term depending on the previous
element of the sequence. The sequence is well defined by using a classical
result of Rozovskii (see Theorem 2, pp. 133, in \cite{Rozovskii}).

The main proof of existence of a truncated solution (see Theorem \ref{section:eusrsw}) builds upon the arguments developed in \cite{CL1} and \cite{CL2},
where the analysis of the vorticity version of the fluid equation was
accomplished. This way the pressure term is eliminated altogether and the
equation was closed by expressing the velocity vector field in terms of the
vorticity via a Biot-Savart operator, (see e.g., \cite{CL1} for
details). This no longer possible here: the equation satisfied by the fluid
vorticity contains a term that depends on $\nabla p$, where $p=\frac{h-b}{%
\epsilon \mathscr{F}}.$ This implies on the one hand that an $L^{p}$%
-transport property for the vorticity is out of reach and, also, that a
control of the higher order derivatives is much more difficult.

As it is well known, compressible systems are much harder to analyse than
their incompressible counterparts. In the particular case of (\ref{shortsrsw}%
), the pressure term is $p=\frac{h-b}{\epsilon \mathscr{F}}$ persists in all
the a priori estimates. To be more precise, the inner product $\langle
v,\nabla p\rangle $ does not vanish as it is the case for incompressible
system. 
Obviously, the same is true higher derivatives $\langle \partial ^{\alpha
}v,\partial ^{\alpha }\nabla p\rangle $. The control of the nonlinear terms
is no longer possible in the same manner as in the incompressible case. As a
result, without adding extra viscosity we cannot "close" any inequalities
involving Sobolev norms of $v$ and $h$: to control $\left\vert \left\vert
v\right\vert \right\vert _{1,2},$ we need a control of $\left\vert
\left\vert h\right\vert \right\vert _{2,2} $ which requires $\left\vert
\left\vert v\right\vert \right\vert _{2,2}~$and so on in a never-ending
procession. In spectral theory language, there is an energy cascade between
low frequency and high frequency modes. To balance the equation and close
the loop, we add viscosity to the system. This way, energy dissipates from
all modes and we succeed to show that, at least for a while, the solution
does not blow-up.

The stochasticity incurs additional technical difficulties. That is because
the solution might blow-up at some random time $\tau _{blow-up}$. For any
deterministic time $t>0$, $\tau _{blow-up}$ might be larger than $t$ for
some realizations of the noise, in other words the solution blows up after
time $t$ if at all, whilst for other realizations of the noise, $\tau
_{blow-up}$ might be smaller than $t$, in other words the solution blows up
before time $t$. As a result, the expectations of the random variables
appearing in (\ref{shortsrsw})\ might not exist. In particular, we cannot
prove that 
\begin{equation}
\mathbb{E}\left[ \left\vert \left\vert a_{t}\right\vert \right\vert
_{1,2}^{2}\right] <\infty ,~~\mathbb{E}\left[ \int_{0}^{t}\left\vert
\left\vert a_{t}\right\vert \right\vert _{2,2}^{2}ds\right] <\infty
,\sum_{i=1}^{\infty }\int_{0}^{t}\mathbb{E}\left[ \left\vert \left\vert 
\mathcal{G}_{i}\left( a_{s}\right) \right\vert \right\vert _{2}\right]
ds<\infty ,  \label{expectations}
\end{equation}%
etc. or any other suitable controls of expectations. In particular, the
Doob-Meyer decomposition of the semi-martingale process $t\rightarrow a_{t}$
contains a \emph{local} martingale. 
On the technical side, the standard (deterministic) Gronwall approach cannot
be used to control the Sobolev norm of the solution $t\rightarrow a_{t}$ of (%
\ref{shortsrsw}) pathwise (because of the stochasticity), neither can a
control of $t\rightarrow \mathbb{E}\left[ \left\vert \left\vert
a_{t}\right\vert \right\vert _{1,2}^{2}\right] $ can be deduced (because of
the possibility of finite blow-up). However, if we "stop" the process at a
time $\tau $ that is \emph{sure to occur} before the blow-up, $\tau <\tau
_{blow-up}$, then we \emph{can} obtain controls on expectations such as the
ones appearing in (\ref{expectations}). 

The definition of the solution of (\ref{shortsrsw}) plays a crucial role in
the analysis, more so than in the deterministic case, and we introduce it in
the next section.

In the absence of noise, one can prove that the viscous RSW has a global
solution for sufficiently small initial datum, $\left\vert \left\vert
a_{0}\right\vert \right\vert \leq \varepsilon $, where $\varepsilon
=\varepsilon \left( \gamma \right) $. That is because one can show that there
exist constants $b=b\left( \gamma \right) $ and $c=c\left( \gamma \right) $ such
that 
\begin{equation*}
{\mathrm{d}}_{\mathrm{t}}\left\vert \left\vert a\right\vert \right\vert
_{1,2}^{2}\leq r_{b,c} \left\vert \left\vert a\right\vert \right\vert
_{1,2}^{2} ,
\end{equation*}%
where $\left\vert \left\vert \cdot \right\vert \right\vert _{1,2}$ is the
Sobolev norm of the space $\left( \mathcal{W}^{1,2}(\mathbb{T}^2)\right) ^{3}$and $%
r_{b,c}:[0,\infty )\rightarrow \lbrack 0,\infty ),\ r_{b,c}\left( s\right)
=bs^{3}-cs$. The reason for this is that the solution of the ODE~$%
d_{t}q=r_{b,c}\left( q\right) $ is an upper bound for $t\rightarrow
\left\vert \left\vert a_{t}\right\vert \right\vert _{1,2}^{2}.$ The function 
$t\rightarrow q_{t}$ is bounded if the initial condition belongs to the
interval $\ [0,\sqrt{c/b}]$. In fact, $\lim_{t\rightarrow \infty }q_{t}=0$
if the initial condition belongs belongs to the interval $\ [0,\sqrt{c/b})$
and it is constant if $q_{0}=\sqrt{c/b}.\ $However, if $q_{0}>\sqrt{c/b}$
then the solution blows up in finite time. This does not necessarily mean
that the solution of the SRSW blows up in finite time. In the stochastic
case, we can deduce a corresponding \emph{stochastic} differential equation
that gives us an upper bound $\left\vert \left\vert a\right\vert \right\vert
_{1,2}^{2}$. The solution of this can be shown to remain bounded only with
positive probability.

Therefore, we can prove that $\mathbb{P}\left( \tau _{blow-up}=\infty
\right) >0$, but not that $\mathbb{P}\left( \tau _{blow-up}=\infty \right) =1$. In
fact, it is not necessarily true that the solution will actually blow up if
it is global:\ in fact we prove that the solution \emph{remains} uniformly
bounded on $[0,\infty )$, not just that it does not blow up in finite time.
The result we obtain gives only a sufficient condition for global existence. 
In future work we aim to show that, under suitable additional assumptions on
the choice of the noise and that of the initial condition the stochastic RSW
equation exists \emph{globally with probability 1}.

Similar results (and similar proofs) hold for $a_{0}\in ( \mathcal{W}^{k,2}(%
\mathbb{T}^2))^{3}, \ k>1$. In this case, the system (\ref{shortsrsw}) has
paths $t\rightarrow a_{t\wedge\tau}$ in the space 
\begin{equation*}
C\left( [0,T],\left( \mathcal{W}^{k,2}(\mathbb{T}^2)\right) ^{3}\right) \cap
L^{2}\left(0,T;\left( \mathcal{W}^{k+1,2}(\T^2)\right) ^{3}\right)
,~\ ~T>0.
\end{equation*}%
In this case the maximal time can be characterized in a similar manner to
the classical Beale-Kato-Majda criterion. The justification of such criteria
is the subject of future work.

We construct sequences of
approximating solutions which converge in a suitable sense to a truncated
form of the original SRSW system. Then we show that the truncation can be
lifted. As opposed to the case of the Euler equation, here we can lift the
truncation only up to a positive stoping time. This is due to a lack of
transport properties for both variables, which derive from the
compressibility condition and the form of the nonlinear terms. Therefore we
obtain a local solution for the original system \eqref{eq:mainsrsw}. This is
proven in Section \ref{section:eusrsw}. In Section \ref%
{subsection:maximalsol} we show that this solution is maximal and due to the
pathwise uniqueness property from Section \ref{sect:pathwiseuniqueness}, it
is also strong in probabilistic sense. In Section \ref%
{section:globaletruncated} we show global existence for the truncated model. In Section \ref{section:analyticalapprox} we study the analytical
properties of the approximating sequence. A couple of a priori estimates and other useful results are presented in the Appendix. 

\section{Preliminaries and notations}
\label{preliminariesSRSW}
In this section we introduce the main notations, together with the It\^{o} form of the system, definition of solutions and other assumptions and remarks. 
\subsection{Notations}
\begin{itemize}
\item In this paper we work on the two-dimensional torus denoted by $\mathbb{T}^{2}
$. 
\item Let $(\Omega ,\mathcal{F},(\mathcal{F}_{t})_{t},\mathbb{P},(W^{i})_{i})$
be a fixed stochastic basis consisting of the filtered probability space $%
(\Omega ,\mathcal{F},(\mathcal{F}_{t})_{t},\mathbb{P})$ and a sequence of
independent one-dimensional Brownian motions $(W^i)_i$ which are adapted to
the complete and right-continuous filtration $(\mathcal{F}_t)_t$. 
\item  Let 
$\mathcal{M}_t:=L^2\left( \Omega, C\left([0,t];\cW^{1,2}(\T^2)^3\right)\right) \cap L^2\left(\Omega, L^2\left(0,t;\cW^{2,2}(\T^2)^3\right)\right)$ for all $t \geq 0$. For a stochastic process $b$ belonging to the space $\mathcal{M}_t$ define the norm
\begin{equation*}\label{t22norm}
    \|b\|_{t,2,2}^2 := \displaystyle\sup_{s\in[0,t]}\|b_s\|_{1,2}^2 + \displaystyle\int_0^t\|b_s\|_{2,2}^2ds.
\end{equation*}
 Define also
\begin{equation*}
    \vertiii{a_t} :=\sup_{s\in 
\left[ 0,t\right] }\left( \Vert v_{s}\Vert _{1,2}^{2}+\Vert h_{s}\Vert
_{1,2}^{2}\right).
\end{equation*}
It follows that $a, b \in \mathcal{M}_{t}$ implies $\mathbb{E}[\|b\|_{t,2,2}] < \infty$ and $\mathbb{E}\left[ \vertiii{a_t}\right] <\infty$ for all $t>0$.
\item Let  $\alpha \in (0,1), p\in [2,\infty)$ and let $H$ be a Hilbert space. Then the fractional Sobolev space $\mathcal{W}^{\alpha,p}(0,T; H)$ is endowed with the norm 
\begin{equation*}
\|f\|_{\mathcal{W}^{\alpha, p}(0, T ; H)}^p:=\int_{0}^{T}\left\|f_{t}\right\|_{H}^{p} d t+\int_{0}^{T} \int_{0}^{T} \frac{\left\|f_{t}-f_{s}\right\|_{H}^{p}}{|t-s|^{1+\alpha p}} d t d s.
\end{equation*} 
\item $C$ is a generic constant and can differ from line to line. 
\end{itemize}

\subsection{It\^{o} form and definition of solutions}
The
expanded version of the stochastic system (\ref{shortsrsw}) is 
\begin{subequations}
\label{eq:mainsrsw}
\begin{equation}
dv_{t}+\big[u_{t}\cdot \nabla v_{t}+f\hat{z}\times u_{t}+\nabla p_{t}\big]dt+%
\displaystyle\sum_{i=1}^{\infty }\big[(\mathcal{L}_{i}+\mathcal{A}_{i})v_{t}%
\big]\circ dW_{t}^{i}=\nu \Delta v_{t}dt  \label{eq:mainsrswa}
\end{equation}%
\begin{equation}
dh_{t}+\nabla \cdot (h_{t}u_{t})dt+\displaystyle\sum_{i=1}^{\infty }\big[%
\nabla \cdot (\xi _{i}h_{t})\big]\circ dW_{t}^{i}=\eta \Delta h_{t}dt.
\label{eq:mainsrswb}
\end{equation}
\noindent

The corresponding It\^{o} form of the system \eqref{eq:mainsrsw}
is given below 
\end{subequations}

\begin{subequations}
\label{eq:itomainsrsw}
\begin{alignat}{2}  \label{eq:itomainsrswb}
& dv_t + \big[\mathcal{L}_{u_t} v_t + f\hat{z} \times u_t + \nabla p_t -
\nu\Delta v_t\big]dt + \displaystyle\sum_{i=1}^{\infty}\big[(\mathcal{L}_i + 
\mathcal{A}_i)v_t\big] dW_t^i = \frac{1}{2}\displaystyle\sum_{i=1}^{\infty}%
\big[(\mathcal{L}_i + \mathcal{A}_i)^2v_t\big]dt &  \\
& dh_t+ \big[\nabla \cdot (h_tu_t) - \eta \Delta h_t\big]dt+\displaystyle%
\sum_{i=1}^{\infty}\mathcal{L}_ih_t dW_t^i=\frac{1}{2}\displaystyle%
\sum_{i=1}^{\infty}\mathcal{L}_i^2h_tdt. & 
\end{alignat}
\end{subequations}

In the following we will work with the It\^{o} version (\ref{eq:itomainsrsw}%
) of system \eqref{eq:mainsrsw}. The definition of a solution of the system (%
\ref{eq:itomainsrsw}) is made explicit in Definition \ref{defsol} below.
Using the It\^{o} version (\ref{eq:itomainsrsw}) of the system as basis for
the well-posedness analysis enables us to match the constraints imposed on
the initial condition to guarantee the well-posedness of the deterministic
system. In particular, Theorems \ref{mainthm} and \ref{truncatedsrsw} below state that the system (%
\ref{eq:itomainsrsw}) is well-posed provided $a_0=(v_0,h_0)\in (W^{1,2}(%
\mathbb{T}^2))^3$.

The required smoothness constraint on $a_0$ that ensures the existence of a
strong solution for the equation \emph{written in Stratonovitch form}, i.e.
the system \eqref{eq:mainsrsw}, is $a_0=(v_0,h_0)\in (W^{1,2}(\mathbb{T}^2))^3$%
. Let us explain why: the Stratonovich integrals in the system %
\eqref{eq:mainsrsw} require that the integrands are semi-martingales, in
this case, the processes $\mathcal{G}_{i} a$. The evolution equation for the
processes $\mathcal{G}_{i} a_s$ involves the terms $\mathcal{G}_{i} ^{3} a_s$
which make sense if the paths of the solution are in $L^2 ([0,T], W^{3,2}(\mathbb{T}^2))$.
This can be achieved, due to the added viscosity term, if the initial
condition $a_0\in W^{2,2}(\mathbb{T}^2)$. To avoid the additional smoothness requirement
we work with the It\^o version (\ref{eq:itomainsrsw}) of the system.

We introduce the following notions of solutions:

\begin{definition}
\label{defsol} $\left.\right.$

\begin{enumerate}
\item[a.] A pathwise \underline{local solution} of the SRSW system is given
by a pair $(a,\tau )$ where $\tau :\Omega \rightarrow \lbrack 0,\infty ]$ is
a strictly positive bounded stopping time, $a_{\cdot \wedge \tau}:\Omega
\times \lbrack 0,\infty )\times \mathbb{T}^{2}\rightarrow \mathbb{R}^{3}$
is an $\mathcal{F}_{t}$-adapted process for any $%
t\geq 0$, with initial condition $a_{0}$,
such that 
\begin{equation*}
a_{\cdot \wedge \tau }\in L^{2}\left( \Omega ;C\left( [0,T];\mathcal{W}%
^{1,2}(\mathbb{T}^{2})^{3}\right) \right) \cap L^{2}\left( \Omega
;L^{2}\left(0,T;\mathcal{W}^{2,2}(\mathbb{T}^{2})^{3}\right) \right)
\end{equation*}
and the SRSW system (\ref{shortsrsw}) is satisfied locally i.e.%
\begin{equation}
a_{t\wedge \tau }=a_{0}+\int_{0}^{_{t\wedge \tau }}\tilde{F}(a_{s})
ds+\sum_{i=1}^{\infty }\int_{0}^{_{t\wedge \tau }}\mathcal{G}_{i}(a_{s})
dW_{s}^{i}=\gamma \int_{0}^{_{t\wedge \tau }}\Delta a_{s}ds,
\end{equation}%
holds $\mathbb{P}$-almost surely, as an identity in $%
L^{2}(\mathbb{T}^{2})^{3},$ with $\tilde{F}(a_s):= F(a_s) + \frac{1}{2}\displaystyle\sum_{i=1}^{\infty}%
\mathcal{G}_i^2(a_s)$. 

\item[b.] If $\tau =\infty $ then the solution $(v,h)$ is called \underline{global}.

\item[c.] A pathwise \underline{maximal solution} of the SRSW system is
given by a pair $(a, \mathcal{T})$ where $\mathcal{T}: \Omega \rightarrow
[0, \infty]$ is a non-negative stopping time and $a=(a_{t\wedge\mathcal{T}})_t$, $a_{t\wedge\mathcal{T}}: \Omega\times \mathbb{T}^2 \rightarrow 
\mathbb{R}^3$ is a process for which there exists an increasing
sequence of stopping times $(\tau^n)_n$ with the following properties: 
\begin{itemize}
\item[i.] $\mathcal{T}=\lim_{n\rightarrow \infty }\tau ^{n}$ and $\mathbb{P}(%
\mathcal{T}>0)=1$

\item[ii.] $(a,\tau ^{n})$ is a pathwise local solution of the SRSW system
for every $n\in \mathbb{N}$

\item[iii.] if $\mathcal{T}<\infty $ then

\begin{equation*}
\displaystyle\limsup_{t \rightarrow \mathcal{T}}  \Vert a_{t}\Vert
_{1,2} = \infty .
\end{equation*}
\end{itemize}

\item[d.] A \underline{weak/martingale local solution} of the SRSW system is
given by a triple \newline
$((\check{a},(\check{W}^{i})_{i}),(\check{\Omega},\mathcal{\check{F}},\check{%
\mathbb{P}}),(\mathcal{\check{F}}_{t})_{t})$ such that $(\check{\Omega},%
\mathcal{\check{F}},(\mathcal{\check{F}}_{t})_{t},\check{\mathbb{P}},(\check{%
W}^{i})_{i})$ is a stochastic basis, $\check{a}$ is a continuous $(\mathcal{%
\check{F}}_{t})_{t}$-adapted real valued process, $\check{a}: \check{\Omega}%
\times \mathbb{T}^{2}\rightarrow \mathbb{R}^{3}$, which satisfies %
\eqref{shortsrsw} for a stopping time $\tau :\check{\Omega}\rightarrow
\lbrack 0,\infty ]$, 
$(\check{W}^{i})_{i}$ are independent $(\mathcal{\check{F}}_{t})_{t}$%
-adapted Brownian motions, and all identities hold $\check{\mathbb{P}}$%
-almost surely in $L^{2}(\mathbb{T}^{2})^{3}$.
\end{enumerate}
\end{definition}

\begin{remark}
\begin{itemize} 

\item We will show below that the SRSW system (\ref{shortsrsw}) satisfies
the local uniqueness property. In other words, if $(a^1,\tau^1)$ and $%
(a^2,\tau^2)$ are two local solutions of system (\ref{shortsrsw}), then they
must coincide on the interval $[0,\tau^1\wedge\tau^2]$. Using the local
uniqueness property, we will deduce, that a stopping time $\mathcal{T}$
satisfying property iii. of the definition above is the largest stopping
time with properties $i.$ and $ii.$, that is for any other pair $%
(a^{\prime},\mathcal{T}^{\prime })$ which satisfies $i.$ and $ii.$ we
necessarily have $\mathcal{T}^{\prime }\leq \mathcal{T}$ $\mathbb{P}$-a.s.,
and $a=a^{\prime }$ on $[0,\mathcal{T}^{\prime}]$. 

\item The first two definitions of solutions are established with respect to
a fixed stochastic basis, the solutions being \textit{strong} in
probabilistic sense. The solution defined at d. is weak in probabilistic
sense, meaning that $(\check{v}, \check{h})$ are not necessarily adapted to
the original filtration $(\mathcal{F}_t)_t$ generated by the driving
Brownian motion which corresponds to the SRSW system \eqref{eq:mainsrsw}. 
\end{itemize}
\end{remark}

\subsection{Assumptions and remarks}

\label{xiassumptionssrsw}

The vector fields $\xi_i:\mathbb{T}^2 \rightarrow \mathbb{R} ^2 $ are chosen
to be time-independent and divergence-free, such that 
\begin{equation}  \label{xiassumptc}
\displaystyle\sum_{i=1}^{\infty} \|\xi_i\|_{4,\infty}^2 < \infty.
\end{equation}
Condition \eqref{xiassumptc} implies
that the infinite sums of stochastic integrals
\begin{equation}  \label{infsum11}
\sum_{i=1}^{\infty} \int_0^t \mathcal{G}_{i} (a_s)dW_s^i, \ \
\sum_{i=1}^{\infty} \int_0^t \nabla \mathcal{G}_{i} (a_s)dW_s^i
\end{equation}
are well defined and belong to $L^2(0, T; L^2(\mathbb{T}^2)^3)$, provided the process $s\rightarrow a_s$ has paths in the space $L^2(0,T; \mathcal{W}^{2,2}(\mathbb{T}^2)^3)$ for $T\geq 0$. Local solutions of the SRSW model as defined above have this property.  
Similarly, the infinite sums of the 
Riemann-Stieltjes integrals
\begin{equation}  \label{infsum111}
\sum_{i=1}^{\infty} \int_0^t \mathcal{G}^2_{i} (a_s)ds, 
\end{equation}
are well-defined and
belong to $L^2(0, T; L^2(\mathbb{T}^2)^3)$. 

\section{Existence and uniqueness of strong pathwise solutions for the SRSW
system}
\label{section:eusrsw} 
In this section we present the main results of this paper. 
\begin{theorem}
\label{mainthm} Let $\mathcal{S} = (\Omega, \mathcal{F}, (\mathcal{F}_t)_t, 
\mathbb{P}, (W^i)_i)$ be a fixed stochastic basis and $a_0\in \mathcal{W}%
^{1,2}(\mathbb{T}^2)^3$. Then the stochastic rotating shallow water system %
\eqref{eq:mainsrsw} admits a unique pathwise maximal
solution $(a,\mathcal{T})$
which belongs to the space $L^{2}\left( \Omega ;C\left( [0,\mathcal{T});\mathcal{W}^{1,2}(\mathbb{T}^{2})^{3}\right) \right) \cap L^{2}\left( \Omega;L^{2}\left(0,\mathcal{T};\mathcal{W}^{2,2}(\mathbb{T}^{2})^{3}\right) \right)$.
\end{theorem}

The existence of a solution for the system \eqref{eq:mainsrsw} is proved by first
showing that a truncated version of it has a solution and then removing the
truncation up to a positive stopping time. In particular, we truncate
the nonlinear terms in \eqref{eq:mainsrsw} using a smooth
function $f_R: \mathbb{R}_{+} \rightarrow [0,1]$ equal to $1$ on $[0,R]$,
equal to $0$ on $[R+1, \infty)$, and decreasing on $[R, R+1]$, 
$f_{R}(a_{t}):=f_{R}\left( \Vert a_t\Vert _{1,2}\right)$
where 
$\Vert a_t\Vert _{1,2}:=\Vert v_t\Vert _{1,2}+\Vert h_{t}\Vert _{1,2}$
for arbitrary $R>0$.
The choice of the truncation $f_{R}(a_{t})$ is such that the nonlinear
terms are uniformly bounded \emph{pathwise} in $L^{2}\left( 0,T;L^{2}(\mathbb{T}%
^{2})^3\right) $ for any $T\geq 0$. Then we have the following: 
\begin{theorem}
\label{truncatedsrsw} Let $\mathcal{S} = (\Omega, \mathcal{F}, (\mathcal{F}%
_t)_t, \mathbb{P}, (W^i)_i)$ be a fixed stochastic basis and $a_0\in \mathcal{W%
}^{1,2}(\mathbb{T}^2)^3$. Then the truncated system 
\begin{subequations}
\label{truncatedsrsweqn}
\begin{alignat}{2}
& dv_{t}^R+\big[f_{R}(a_{t}^R)\mathcal{L}_{u_{t}^R}v_{t}^R+f\hat{z}\times
u_{t}^R+ \nabla p_{t}^R+\nu \Delta v_{t}^R\big]dt+\displaystyle%
\sum_{i=1}^{\infty }\big[(\mathcal{L}_{i}+\mathcal{A}_{i})v_{t}^R\big]\circ
dW_{t}^{i}=0 & &  \label{truncatedsrsweqnh} \\
& dh_{t}^R+\big[f_{R}(a_{t}^R)\nabla \cdot (h_{t}^Ru_{t}^R)+\eta \Delta h_{t}^R\big]dt+%
\displaystyle\sum_{i=1}^{\infty }\big[\nabla \cdot (\xi _{i}h_{t}^R)\big]\circ
dW_{t}^{i}=0 & &
\end{alignat}%
admits a unique global pathwise solution $a^{R}=\{(v_{t}^{R},h_{t}^{R}),t\in
\lbrack 0,\infty )\}$ such that 
\end{subequations}
\begin{equation*}
a^{R}\in L^{2}\left( \Omega ;C\left( [0,T];\mathcal{W}^{1,2}(\mathbb{T}%
^{2})^{3}\right) \right) \cap L^{2}\left( \Omega ;L^{2}\left( 0,T;\mathcal{W%
}^{2,2}(\mathbb{T}^{2})^{3}\right) \right)
\end{equation*}%
for any $T>0$.
\end{theorem}
Theorem \ref{truncatedsrsw} is proved in Section \ref{section:globaletruncated}. Define
\begin{equation}\label{taur}
\tau ^{R}:=\inf_{t\geq 0}\left\{\Vert a_t\Vert _{1,2}\geq R\right\}.
\end{equation}

\begin{proposition}
\label{prop:localtoglobal} 
Given $a_0 \in \left( \cW^{1,2}(\T^2)\right)^3$ and $R>0$, the restriction $a: \Omega \times [0, \tau_R) \times \T^2 \rightarrow \mathbb{R}^3$ of the global solution $a : \Omega \times [0, \infty) \times \T^2 \rightarrow \mathbb{R}^3$ corresponding to the truncated system \eqref{truncatedsrsweqn} is a local solution of the original SRSW system \eqref{eq:mainsrsw}. 
\end{proposition}

\noindent\textbf{\textit{Proof}} For $t\in[0, \tau_R]$ $f_R(a^R)=1$ therefore the truncated system \eqref{truncatedsrsweqn} and the original
SRSW system \eqref{eq:mainsrsw} coincide.

\subsection{Maximal solution for the SRSW system}
\label{subsection:maximalsol}

\begin{proposition}\label{prop:max.sol.}
Given $a_0\in \mathcal{W}^{1,2}(\mathbb{T}^2)^3$ and $R>0$, there exists a
unique maximal solution $\left( a,\mathcal{T}\right) $ of the original SRSW
system \eqref{eq:mainsrsw} such that 
\begin{equation}
\displaystyle\limsup_{t\rightarrow \mathcal{T}}\left\vert \left\vert
a_t\right\vert \right\vert _{1,2}=\infty,  \label{limsup}
\end{equation}%
whenever $\mathcal{T}<\infty $.
\end{proposition}

\noindent \textbf{\textit{Proof.}}\\
\textbf{Existence.} If we choose $R=n$ in Proposition \ref%
{prop:localtoglobal} then $\left( a^{n},\tau ^{n}\right) $ is a local
solution of the SRSW system \eqref{eq:mainsrsw}. Moreover, observe that $%
a^{n+1}$ satisfies equation \eqref{truncatedsrsweqn} for $R=n$ on the
interval $[0,\tau ^{n}]$. By the local uniqueness, it follows that 
\begin{equation}
a^{n+1}|_{[0,\tau ^{n}]}=a^{n}|_{[0,\tau ^{n}]}.  \label{consistent}
\end{equation}%
Define $\mathcal{T}:=\displaystyle\lim_{n\rightarrow \infty }\tau ^{n}$ 
and 
\begin{equation}
a|_{[0,\tau ^{n}]}:=a^{n}|_{[0,\tau ^{n}]}.  \label{maximalsolution}
\end{equation}%
Definition (\ref{maximalsolution}) is consistent following (\ref%
{consistent}). It only remains to show (\ref{limsup}). 
If $\displaystyle\lim_{n\rightarrow \infty }\tau ^{n}=\mathcal{T}<\infty $ 
then
\begin{equation*}
\displaystyle\limsup_{t\rightarrow \mathcal{T}}\left\vert \left\vert
a_t\right\vert \right\vert _{1,2}\ge    \displaystyle\limsup_{n\rightarrow\infty}\|a_{\tau_n}\|_{1,2} =  \displaystyle\limsup_{n\rightarrow\infty} n = \infty.
\end{equation*}
\textbf{Uniqueness}. 
Assume that $\left( \bar{a}{,\mathcal{\bar{T}}}\right) $ is another solution with 
$\left( \bar{a},\bar\tau_n \right) $, $n\ge 1$ being the corresponding sequence of local solutions converging to the maximal solution. 
By the uniqueness of the truncated
equation it follows that $\bar{a}=a$ on $[0,\bar{\tau}^{n}\wedge \tau ^{n}]$%
. By taking the limit as $n\rightarrow \infty $ it follows that $\bar{a}=a$
on $[0,\mathcal{T}\wedge \mathcal{\bar{T}})$. We prove next that $\mathcal{T}%
=\mathcal{\bar{T}}$, $\mathbb{P}-a.s.$. Let us assume that 
\begin{equation*}
\mathbb{P}\left( \omega \in \Omega ,\mathcal{T}\left( \omega \right) \neq 
\mathcal{\bar{T}}\left( \omega \right) \right) >0.
\end{equation*}%
Observe that 
\begin{equation*}
\mathbb{P}\left( \omega \in \Omega ,\mathcal{T}\left( \omega \right) \neq 
\mathcal{\bar{T}}\left( \omega \right) \right) \leq \mathbb{P}\left( \Xi
^{1}\right) +\mathbb{P}\left( \Xi ^{2}\right) ,
\end{equation*}%
where 
\begin{eqnarray*}
\Xi ^{1} &=&\left\{ \omega \in \Omega ,\mathcal{T}\left( \omega \right)
<\infty ,\mathcal{T}\left( \omega \right) <\mathcal{\bar{T}}\left( \omega
\right) \right\} , \\
\Xi ^{2} &=&\left\{ \omega \in \Omega ,\mathcal{\bar{T}}\left( \omega
\right) <\infty ,\mathcal{\bar{T}}\left( \omega \right) <\mathcal{T}\left(
\omega \right) \right\} .
\end{eqnarray*}%
We prove that $\mathbb{P}\left( \Xi ^{1}\right) =\mathbb{P}\left( \Xi
^{2}\right) =0$. The two sets are symmetric so we show this only for the
first one. From the definition of the local solution, observe that%
\begin{equation*}
\mathbb{E}\left[ \sup_{s\in \lbrack 0,\bar{\tau}^{n}\left( \omega \right)
]}\Vert \bar{a}_t\Vert_{1,2}\right] <\infty 
\end{equation*}%
hence%
\begin{equation*}
\mathbb{P}\left( \omega \in \Omega ,\sup_{s\in \lbrack 0,\bar{\tau}%
^{n}\left( \omega \right) ]}\Vert \bar{a}_t\Vert_{1,2}=\infty \right) =0.
\end{equation*}%
However, if $\mathcal{T}\left( \omega \right) <\infty $ and $\mathcal{T}%
\left( \omega \right) <\bar{\tau}^{n}\left( \omega \right) $ then 
\begin{equation*}
\infty =\sup_{s\in \lbrack 0,\mathcal{T}\left( \omega \right) )}\Vert a_s\Vert
_{1,2}=\sup_{s\in \lbrack 0,\mathcal{T}\left( \omega \right) )}\Vert \bar{a_s
}\Vert_{1,2}\leq \sup_{s\in \lbrack 0,\bar{\tau}^{n}\left( \omega \right)
]}\Vert \bar{a}_s\Vert_{1,2}.
\end{equation*}%
It follows that $\mathbb{P}\left( \omega \in \Omega ,\mathcal{T}\left(
\omega \right) <\infty ,\mathcal{T}\left( \omega \right) <\bar{\tau}%
^{n}\left( \omega \right) \right) =0.$ Hence 
\begin{equation*}
\mathbb{P}\left( \Xi ^{1}\right) =\lim_{n\rightarrow \infty }\mathbb{P}%
\left( \omega \in \Omega ,\mathcal{T}\left( \omega \right) <\infty ,\mathcal{%
T}\left( \omega \right) <\bar{\tau}^{n}\left( \omega \right) \right) =0.
\end{equation*}%
This completes the proof of the uniqueness claim. \\

\noindent The purpose of the next proposition is to show that the maximal solution constructed in Proposition \ref{prop:max.sol.} has paths in $L^{2,loc}\left( [0,\cT), (\mathcal{W}^{2,2}(\mathbb{T}^2))^3\right)$. Recall the definition of $\tau^N$ as given in \eqref{taur} and the definition of the $\|a\|_{t,2,2}$ as introduced in \eqref{t22norm} and introduce a new sequence of stopping times 
\begin{equation*}
    \hat{\tau}^M := \displaystyle\inf_{t\geq 0}\{ \|a\|_{t,2,2} \geq M\}.
\end{equation*}
Define 
$\hat{\cT} := \displaystyle\lim_{M\rightarrow\infty} \hat{\tau}^M$.

\begin{corollary}\label{stoppedestim}
Let $(a, \cT)$ be the maximal solution constructed in Proposition \ref{prop:max.sol.}. Then the process  $t\rightarrow a_{t\wedge \tau^R}$ takes values in 
\begin{equation*}
L^{2}\left( \Omega ;C\left( [0,T];\mathcal{W}^{1,2}(\mathbb{T}%
^{2})^{3}\right) \right) \cap L^{2}\left( \Omega ;L^{2}\left( 0,T;\mathcal{W%
}^{2,2}(\mathbb{T}^{2})^{3}\right) \right)
\end{equation*}%
for any $T>0$. In particular, 
\begin{equation}\label{difnorm}
\mathbb{E}[\|a\|_{\tau^R\wedge T,2,2}^2] <\infty
\end{equation}
for any $R,T>0$.
\end{corollary}
\noindent\textbf{Proof.} Immediate from Theorem  \ref{truncatedsrsw} after observing that $a=a^R$ on $[0,\tau^R]$ and $a^R \in \mathcal{M}_T$.\\[2mm]

\noindent We are now ready to show the equality between the two stopping times $\cT$ and $\bar\cT$.
\begin{proposition}\label{highernormcontrol}
Let $(a, \cT)$ be the maximal solution constructed in Proposition \ref{prop:max.sol.}. Then 
\begin{equation*}
    \mathbb{P}\left( \cT = \hat{\cT}  \right) = 1. 
\end{equation*}
\end{proposition}
\noindent\textbf{Proof.}
 One can observe that $\hat{\tau}^M \leq \tau^M$ for all $M\geq 0$, and therefore $\hat{\cT} \leq \cT$ $\mathbb{P}-a.s.$. We show that $\cT \leq \hat{\cT}$ $\mathbb{P}-a.s.$. On the set $\{ \hat{\cT} = \infty\}$ the inequality is trivially true, so we only need to show that
 \begin{equation*}
     \mathbb{P}\left( \{\omega \in \Omega: \hat{\cT}(\omega) < \infty, \cT(\omega) \leq \hat{\cT}(\omega)\} \right) =1. 
 \end{equation*}
 Note that
 \begin{equation*}
     \begin{aligned}
    \{ \omega \in \Omega: \displaystyle \|a(\omega)\|_{\tau^N,2,2} < \infty \} = \displaystyle\bigcup_M  \{ \omega \in \Omega:  \|a(\omega)\|_{\tau^N,2,2} < M \} = \displaystyle\bigcup_M \{ \omega \in \Omega: \tau^N(\omega) <  \hat{\tau}^M(\omega) \}
     \end{aligned}
 \end{equation*}
 and
 \begin{equation*}
     \begin{aligned}
       \displaystyle\bigcup_M \{ \omega \in \Omega: \tau^N(\omega) < \hat{\tau}^M(\omega)\} \subset \{ \omega \in \Omega: \hat{\cT}(\omega)<\infty,\tau^N(\omega) < \hat{\cT}(\omega) \}.
     \end{aligned}
 \end{equation*}
From Corollary \ref{stoppedestim} we deduce that 
\begin{equation*}
    \begin{aligned}
    \mathbb{P}\left( \{ \omega \in \Omega:  \|a(\omega)\|_{\tau^N\wedge T,2,2} < \infty\} \right) =1, \ \ \ \forall N\in\mathbb{N}, \ \ \ T>0.
    \end{aligned}
\end{equation*}
and, since
$\limsup_{t\rightarrow \hat{\cT}} \|a(\omega)\|_{t,2,2}=\infty$ on the set $\{\hat{\cT}<\infty\}$, we deduce that 
\[
    \mathbb{P}\left( \{\omega \in \Omega: \hat{\cT}(\omega)<\infty,  \tau^N(\omega) \wedge T< \hat{\cT}(\omega)\}\right)=1,
    \ \ \ \forall N\in\mathbb{N}, \ \ \ T>0,
\]
therefore
\begin{equation*}
    \begin{aligned}
    \mathbb{P}\left( \{\omega \in \Omega: \hat{\cT}(\omega)<\infty,  \tau^N(\omega) <\hat{\cT}(\omega)\}\right)=
    \mathbb{P}\left( \bigcap_{L>1}\{\omega \in \Omega: \hat{\cT}(\omega)<\infty,  \tau^N(\omega)\wedge L < \hat{\cT}(\omega)\}\right) = 1. 
    \end{aligned}
\end{equation*}
Then we have
\begin{equation*}
    \begin{aligned}
    \mathbb{P}\left( \{\omega \in \Omega: \hat{\cT}(\omega)<\infty, \cT(\omega) \leq \hat{\cT}(\omega)\}\right) &\ge \mathbb{P}\left( \displaystyle\bigcap_N \{ \omega \in \Omega: \hat{\cT}(\omega)<\infty, \tau^N(\omega) < \hat{\cT}(\omega)\} \right) 
     =1. 
    \end{aligned}
\end{equation*}

\subsection{Pathwise uniqueness for the truncated SRSW system}
\label{sect:pathwiseuniqueness}

 Let $a^{R,1}=(v^{R,1}, h^{R,1})$ and $a^{R,2}=(v^{R,2}, h^{R,2})$ be two solutions of the truncated system starting from the initial conditions $a_0^1$, $a_0^2 \in \mathcal{W}^{1,2}(\mathbb{T}^2)^3$, respectively. We denote the corresponding differences by $\bar{a}^R := a^{R,1} - a^{R,2}$, $\bar{v}^R:= v^{R,1}-v^{R,2}$, $%
\bar{h}^R:=h^{R,1}-h^{R,2}$. Also $\bar{u}^R:= u^{R,1}-u^{R,2}$, $\bar{p}^R:= p^{R,1}-p^{R,2}$. Assume that $\tau_M^{R,i}$ are the stopping times defined as \begin{equation*}
    \tau_M^{R,i}:= \displaystyle\inf_t\{ t\geq 0, \|a^{R,i}\|_{t,2,2}\geq M \}.
\end{equation*}
Define $\bar{\tau}_M^R := \tau_M^{R,1}\wedge \tau_M^{R,2}$. 

\begin{remark}\label{truncsol}
We have $\displaystyle\lim_{M\rightarrow\infty}\tau_M^{R,i}=\infty \ \ \mathbb{P}-a.s.$. This is because
\begin{equation*}
    \mathbb{P}\left(\tau_M^{R,i}\leq N\right) = \mathbb{P}(\|a^{R,i}\|_{N,2,2} \geq M) \leq \frac{\mathbb{E}\left[\|a^{R,i}\|_{N,2,2}^2\right]}{M^2}
\end{equation*}
hence
\begin{equation*}
    \mathbb{P}\left(\displaystyle\lim_{M\rightarrow\infty}\tau_M^{R,i}\leq N \right) \leq \displaystyle\lim_{M\rightarrow\infty}\mathbb{P}\left(\tau_M^{R,i}\leq N\right) = 0. 
\end{equation*}
Then
\begin{equation*}
    \begin{aligned}
    \mathbb{P}\left( \displaystyle\lim_{M\rightarrow\infty}\tau_M^{R,i}=\infty\right) & = 1 - \mathbb{P}\left( \displaystyle\lim_{M\rightarrow\infty}\tau_M^{R,i} < \infty\right)
    \geq 1 - \displaystyle\sum_N\mathbb{P}\left( \displaystyle\lim_{M\rightarrow\infty}\tau_M^{R,i} < N\right)
    = 1. 
    \end{aligned}
\end{equation*}
Consequently, also $\displaystyle\lim_{M\rightarrow\infty}\bar{\tau}_M^R = \infty$. 
\end{remark}
\begin{theorem}\label{prop:uniquenesstruncated} 
Let $a^{R,1}$, $a^{R,2}$ be two solutions of the truncated SRSW system 
\eqref{truncatedsrsweqn},
which take values in the space 
$L^2\left( \Omega, C\left([0,T],\cW^{1,2}(\T^2)^3\right)\right) \cap L^2\left(\Omega, L^2\left(0,T;\cW^{2,2}(\T^2)^3\right)\right)$
and start from the initial conditions $a_0^1$, $a_0^2 \in \mathcal{W}^{1,2}(\mathbb{T}^2)^3$, respectively. 
Then there exists a constant $C=C(M)$ 
such that 
\begin{equation*}
    \mathbb{E}\left[\|\bar{a}_{t\wedge\bar{\tau}_M^R}^R\|_{1,2}^2\right] \leq Ce^{Ct}\|\bar{a}_0\|_{1,2}^2,
\end{equation*}
where $\bar{a}^R := a^{R,1} - a^{R,2}$ and $\bar{\tau}_M^R := \tau_M^{R,1}\wedge \tau_M^{R,2}.$
In particular, following from Remark \ref{truncsol}, the truncated SRSW system \eqref{truncatedsrsweqn} has a unique solution in the space
\begin{equation*}
L^2\left( \Omega, C([0,T],\cW^{1,2}(\T^2))\right) \cap L^2\left(\Omega, L^2(0,T;\cW^{2,2}(\T^2))\right).
\end{equation*}
\end{theorem}

\noindent\textbf{\textit{Proof}} 
\noindent We show that
\begin{equation}\label{unique1}
    \begin{aligned}
    d\|\bar{a}_t^R\|_{1,2}^2 \leq C(\epsilon,R)\|Z_t\| \|\bar{a}_t^R\|_{1,2}^2 dt + dB_t
    \end{aligned}
\end{equation}
where $\epsilon > 0$,
\begin{equation*}
    \|Z_t\|:= C\left(\|a_t^{R,1}\|_{1,2}^4 + \|a_t^{R,2}\|_{1,2}^4\right),
\end{equation*}
and $dB_t$ is a local martingale given by 
\begin{equation}
    dB_t := 2\displaystyle\sum_{i=1}^{\infty}\left( \langle \bar{v}_t^R, \mathcal{G}_i\bar{v}_t^R\rangle +  \langle \bar{h}_t^R, \mathcal{L}_i \bar{h}_t^R\rangle 
 + \langle \Delta \bar{v}_t^R, \mathcal{G}_i\bar{v}_t^R\rangle  +  \langle \Delta \bar{h}_t^R, \mathcal{L}_i\bar{h}_t^R \rangle \right)
dW_t^i.
\end{equation}
Then
\begin{equation*}
    \begin{aligned}
    \mathbb{E}\left[ e^{-C\displaystyle\int_0^{t\wedge\bar{\tau}_M^R}\|Z_s\|ds} \|\bar{a}_{t\wedge\bar{\tau}_M^R}^R\|_{1,2}^2\right] 
    & \leq \|\bar{a}_0\|_{1,2}^2 + \mathbb{E}\left[\displaystyle\int_0^{t\wedge\bar{\tau}_M^R}e^{-C\displaystyle\int_0^{s\wedge\bar{\tau}_M^R}\|Z_r\|dr}dB_{s}\right]
    \end{aligned}
\end{equation*}
that is

\begin{equation*}
    \begin{aligned}
   \mathbb{E}\left[\|\bar{a}_{t\wedge\bar{\tau}_M^R}^R\|_{1,2}^2\right] & 
    \leq e^{CM^4t}\|\bar{a}_0\|_{1,2}^2
    \end{aligned}
\end{equation*}
since the stopped process $B_{t\wedge\bar{\tau}_M^R}$ is a martingale. By choosing two solutions of the truncated SRSW system \eqref{truncatedsrsweqn} which start from the same initial conditions, we deduce that $\mathbb{P}(a_s^{R,1} = a_s^{R,2}, \forall s\in [0, \bar{\tau}_M^R])=1$ for any $M>0$, that is the two solutions are indistinguishable. Since $\displaystyle\lim_{M\rightarrow\infty}\bar{\tau}_M^R = \infty$ we deduce that the solution is unique globally.
We will now prove that \eqref{unique1} holds, using Lemma \ref{abstractlemmauniqueness}. 
\noindent We can write
\begin{subequations}
\begin{alignat}{2}  
& d\bar{v}_t^R = \left(Q_{\bar{v}^R} - fk\times \bar{v}_t^R + \nu\Delta\bar{v}_t^R + g\nabla%
\bar{p}_t^R - \frac{1}{2}\displaystyle\sum_{i=1}^{\infty}(\mathcal{L}_i+\mathcal{A}_i)^2\bar{v}_t^R \right)dt + \displaystyle\sum_{i=1}^{\infty}(\mathcal{L}_i+\mathcal{A}_i)\bar{v}_t^RdW_t^i  \\
&  d\bar{h}_t^R = \left( Q_{\bar{h}^R} + \eta\Delta \bar{h}_t^R - \frac{1%
}{2}\displaystyle\sum_{i=1}^{\infty}\mathcal{L}_i^2\bar{h}_t^R \right)dt +  \displaystyle\sum_{i=1}^{\infty}\mathcal{L}_i\bar{h}_t^R dW_t^i . 
\end{alignat}
\end{subequations}
where 
$$Q_{\bar{v}^R}:=f_R(a_R^1)u_t^{R,1}\cdot \nabla v_t^{R,1} -
f_R(a_R^2)v_t^{R,2}\cdot \nabla u_t^{R,2} $$
$$Q_{\bar{h}^R} := f_R(a_R^1)\nabla
\cdot (h_t^{R,1}u_t^{R,1}) - f_R(a_R^2)\nabla \cdot
(h_t^{R,2}u_t^{R,2}).$$

\noindent By the It\^{o} formula 
\begin{equation*}
\begin{aligned}
d\|\bar{a}_t^R\|_{1,2}^2  &+ 2\gamma\|\bar{a}%
_t^R\|_{2,2}^2dt  \leq 2\left(\langle \bar{v}_t^R +  \Delta \bar{v}_t^R, Q_{\bar{v}_t^R}\rangle + \langle \bar{h}_t^R +  \Delta \bar{h}_t^R, Q_{\bar{h}_t^R}\rangle- \langle \bar{v}_t^R + \Delta \bar{v}_t^R ,fk\times\bar{v}_t^R + g\nabla\bar{p}_t^R\rangle \right)dt \\
 &+ \displaystyle\sum_{i=1}^{\infty} \left(\langle (\mathcal{L}_i+\mathcal{A}_i)\bar{v}_t^R, (\mathcal{L}_i+\mathcal{A}_i) \bar{v}_t^R \rangle  +  \langle \mathcal{L}_i\bar{h}_t^R, \mathcal{L%
}_i\bar{h}_t^R\rangle  +  \langle \Delta \bar{v}_t^R,(\mathcal{L}_i+\mathcal{A}_i)^2\bar{v}_t^R\rangle  +  \langle \Delta \bar{h}_t^R, \mathcal{L}_i^2\bar{h}%
_t^R\rangle \right)dt \\
& + 2\displaystyle\sum_{i=1}^{\infty}\left( \langle \bar{v}_t^R, (\mathcal{L}_i+\mathcal{A}_i) \bar{v}_t^R\rangle +  \langle \bar{h}_t^R, \mathcal{L}_i \bar{h}_t^R\rangle + \langle \Delta \bar{v}_t^R, (\mathcal{L}_i+\mathcal{A}_i)\bar{v}_t^R\rangle +  \langle \Delta \bar{h}_t^R, \mathcal{L}_i\bar{h}_t^R \rangle \right)
dW_t^i.  
\end{aligned}
\end{equation*}
All the terms which do not contain a stochastic integral are controlled as functions of $C(\zeta,R)\|Z\|\|\bar{a}\|_{1,2}^2 + \zeta \|\bar{a}\|_{2,2}^2$ using Lemma \ref{abstractlemmauniqueness} and Lemma \ref{abstractlemmaexistence}. We choose $\zeta < \gamma $ such that all the terms which are controlled by $\zeta\|\bar{a}\|_{2,2}^2$ on the right hand side cancel out the term $2\gamma\|\bar{a}\|_{2,2}^2$ on the left hand side. Then \eqref{unique1} holds as requested and therefore the two solutions are indistinguishable as processes with paths in $L^2\left( \Omega, C\left([0,T],\cW^{1,2}(\T^2)^3\right)\right) \cap L^2\left(\Omega, L^2\left(0,T;\cW^{2,2}(\T^2)^3\right)\right)$. 

\begin{remark}\label{globsto}
From Proposition \ref{highernormcontrol}, we deduce that $\displaystyle\lim_{M\rightarrow\infty}\tau_M^{i}=\tilde{\cT}^i\ \ \mathbb{P}-a.s.$, for $i=1,2$. 
Consequently, also $\displaystyle\lim_{M\rightarrow\infty}\bar{\tau}_M = \tilde{\cT}^1 \wedge \tilde{\cT}^2$. Moreover, $a^i=a^{M,i}$ on $[0,\tau^{M,i}]$ for $i=1,2$ and arbitrary $M>0$, therefore $\tau^i_M=\tau^{M,i}_M$ and $\bar{\tau}_M := \bar{\tau}_M^2$.
\end{remark}

\begin{corollary}
Let $(a^1, \mathcal{T}^1)$ and $(a^2, \mathcal{T}^2)$ be two maximal solutions of the original system, starting from $a_0^1, a_0^2 \in \mathcal{W}^{1,2}(\mathbb{T}^2)^3$, respectively. Then there is a constant $C=C(M)$ such that 
\begin{equation*}
    \mathbb{E}\left[\|\bar{a}_{t\wedge\bar{\tau}_M}\|_{1,2}^2\right] \leq Ce^{Ct}\|a_0^1-a_0^2\|_{1,2}^2.
\end{equation*}
\end{corollary}
\noindent\textbf{\textit{Proof}} From Remark \ref{globsto} and Theorem \ref{prop:uniquenesstruncated} we deduce that
\begin{equation*}
    \mathbb{E}\left[ \|\bar{a}_{t\wedge\bar{\tau}_M}\|_{1,2}^2\right] = \mathbb{E}\left[\|\bar{a}_{t\wedge\bar{\tau}_M^M}\|_{1,2}^2\right] \leq Ce^{Ct}\|a_0^1-a_0^2\|_{1,2}^2.
\end{equation*}
\begin{remark}
Note that $\displaystyle\lim_{M\rightarrow\infty}\tau_M^i = \tau^i$ (the maximal time of existence) so the continuity covers the common interval of existence.
\end{remark}

\subsection{Global existence for the truncated SRSW system} \label{section:globaletruncated}
\begin{proposition} \label{prop:existencetruncated}
The truncated SRSW system \eqref{truncatedsrsweqn} admits a global solution $a^R=(v^R, h^R)$ such that $a^R_{[0,T]} \in \mathcal{M}_T$ for any $T\geq 0$. In other words
\begin{equation*}
    a^R_{[0,T]} \in L^2\left( \Omega, C([0,T],\cW^{1,2}(\T^2))\right) \cap L^2\left(\Omega, L^2(0,T;\cW^{2,2}(\T^2)^3)\right)
\end{equation*}
for any $T>0$. Moreover
\begin{equation*}
    a^R_{[0,T]} \in L^p\left( \Omega, \mathcal{W}^{\alpha,p}([0,T],L^{2}(\T^2)^3)\right)
\end{equation*}
for any $p\in (2, \infty)$ and $\alpha \in \left[0, \frac{1}{2}\right )$ such that $p\alpha > 1$ 
and 
\begin{equation*}
    a^R_{[0,T]} \in L^p\left( \Omega, C([0,T],\cW^{1,2}(\T^2))\right)
\end{equation*}
for any $T>0$. 
\end{proposition}

\noindent  In the following we will omit the dependence of the truncated system $a^R=(v^R,h^R)$ on $R$ and simply use the notation $a=(v,h)$ to denote it.
The strategy for proving that the truncated system %
\eqref{truncatedsrsw} has a solution is to construct an approximating system
of processes that will converge in distribution to a solution of %
\eqref{truncatedsrsw}. This justifies the existence of a weak solution.
Together with the pathwise uniqueness of the solution of this equation, we
then deduce that strong/pathwise existence holds. \newline

\noindent Recall that $(v_{0},h_{0})\in \mathcal{W}^{1,2}(\mathbb{T}%
^{2})^2\times \mathcal{W}^{1,2}(\mathbb{T}^{2})$. We construct the sequence $%
(v^{n},h^{n})_{n\geq 0}$ with $v_{t}^{0}:=u_{0}^{0},$ $h_{t}^{0}:=h_{0}^{0}$%
, and for $n\geq 1,$ we define $(v^{n},h^{n})_{n\geq 0}$ \ as the solution of
the \textit{linear} SPDE 
\begin{equation*}
\begin{aligned}
&dv_{t}^{n}=\nu \Delta
v_{t}^{n}dt+P_{t}^{n-1,n}(v_{t}^{n})dt-\sum_{i=1}^{\infty }(\mathcal{L}_{i}+%
\mathcal{A}_{i})v_{t}^{n}dW_{t}^{i,n}  \label{itsystitoSRSW} \\
&dh_{t}^{n}=\delta \Delta
h_{t}^{n}dt+Q_{t}^{n-1,n}(h_{t}^{n})dt-\sum_{i=1}^{\infty }\nabla \cdot (\xi
_{i}h_{t}^{n})dW_{t}^{i,n},
\end{aligned}
\end{equation*}%
 where $%
P_{t}^{n-1,n}(v_{t}^{n})$ and $Q_{t}^{n-1,n}(v_{t}^{n})$ are defined,
respectively, as follows (for $t\geq 0$): 
\begin{equation*}
\begin{aligned}
P_{t}^{n-1,n}(v_{t}^{n}) & :=-f_{R}(a_{t}^{n-1})\mathcal{L}%
_{u_{t}^{n-1}}v_{t}^{n-1}-f\hat{z}\times u_{t}^{n}-\nabla p_{t}^{n}+\frac{1}{%
2}\sum_{i=1}^{\infty }(\mathcal{L}_{i}+\mathcal{A}_{i})^{2}v_{t}^{n} \\
& Q_{t}^{n-1,n}(h_{t}^{n}) :=-f_{R}(a_{t}^{n-1})(\nabla \cdot
\left( h_{t}^{n-1}u_{t}^{n-1}\right) )+\frac{1}{2}\sum_{i=1}^{\infty }%
\mathcal{L}_{i}^{2}h_{t}^{n}
\end{aligned} 
\end{equation*}

\begin{theorem}\label{psnsSRSW}
The approximating system admits a unique global solution in the
space 
\begin{equation*}
(v^{n},h^{n})\in L^{2}\left( \Omega ;C\left( [0,T];\mathcal{W}^{1,2}(%
\mathbb{T}^{2})^{3}\right) \right) \cap L^{2}\left( \Omega ;L^{2}\left( 0,T;%
\mathcal{W}^{2,2}(\mathbb{T}^{2})^{3}\right) \right)
\end{equation*}%
and for any $p\geq 2$ there exists a constant $\mathcal{B}_{3}(T,R)$ independent of $n$ such
that 
\begin{equation}
\mathbb{E}\left[ \| (v^{n},h^{n})\|_{T,2,2}^{p}\right] \leq \mathcal{B}_{3}(T,R).  \label{c1}
\end{equation}
 Moreover $%
(v^{n},h^{n})\in L^{p}\left( \Omega ;\mathcal{W}^{\alpha,p}\left( [0,T],
L^{2}(\mathbb{T}^{2})^{3}\right) \right) $ with $p \in(2,\infty),\alpha \in[0,\frac{1}{2})$ such
that $p\alpha \geq 1$and there exists a constant $\mathcal{B}_{4}(p,\alpha
,T,R)$ independent of $n$ such that%
\begin{equation}
\mathbb{E}\left[ \| (v^{n},h^{n})\|_{\mathcal{W}^{\alpha, p}\left( [0,T],L^{2}(\mathbb{T}^{2})^{3}\right) }^p%
\right] \leq \mathcal{B}_{4}(p,T, R).  \label{c2}
\end{equation}
\end{theorem}

\noindent The proof of this theorem is provided in Section \ref{section:analyticalapprox} below.
\begin{proposition}
\label{tightnessISRSW} The family of probability distributions of the solutions $%
(v^{n},h^{n})_{n} $ is relatively compact in the space of probability
measures over $L^{2}\left( \Omega ;C\left( [0,T];L^{2}(\mathbb{T}%
^{2})^{3}\right) \right) $ for any $T\geq 0$.
\end{proposition}

\noindent Proposition \ref{tightnessISRSW} is proven in Section \ref{section:analyticalapprox}. \\

\noindent\textbf{\textit{Proof of Proposition \ref{prop:existencetruncated}.}}
It is in the proof of this proposition that we see the additional difficulties encountered for stochastic models as compared to the 
deterministic models. Let us explain why this is the case. Recall that  Proposition \ref{tightnessISRSW} tells us that the family of 
\emph{probability distributions} of the  approximate solutions $%
(v^n,h^n)_n$ is relatively compact over
$L^{2}\left( \Omega ;C\left( [0,T];L^{2}(\mathbb{T}%
^{2})^{3}\right) \right) $ for any $T\geq 0$. This does not mean that 
the processes themselves are relatively compact. Therefore, in contrast to the deterministic case, we cannot extract a subsequence from $(v^n,h^n)_n$ that will converge pathwise. We can only extract a subsequence $(v^{n_j},h^{n_j})$ that will  converge in distribution over $L^{2}\left( \Omega ;C\left( [0,T];L^{2}(\mathbb{T}^{2})^{3}\right) \right) $ for any $T\geq 0$. We can then find a \emph{different} probability space 
$(\tilde\Omega,\tilde{\mathcal F}, \tilde{\mathbb{P}})$ 
on which we can build copies of $(v^{n_j},h^{n_j})$ with the same distributions as the original ones that converge in 
$L^{2}\left( \tilde\Omega ;C\left( [0,T];L^{2}(\mathbb{T}%
^{2})^{3}\right) \right) $ and, possibly by using a further subsequence, we can also assume that the convergence is pathwise.  This is done by means 
of a classical probabilistic result called the Skorokhod representation theorem, see for example \cite{Billingsley} Section 6, pp. 70. 

Further complications need to be sorted: It is not enough to transfer just the processes $(v^{n_j},h^{n_j})$. The driving Brownian motions $(W_i)_{i=1}^{\infty}$
will need to be mirrored in the new space $(\tilde\Omega,\tilde{\mathcal F}, \tilde{\mathbb{P}})$ as the "mirroring processes" is done for each individual term of the subsequence. We end up with a set of Brownian motions that are different for each element of the sequence, even if we start with a subsequence that is driven by the same set of Brownian 
motions (therefore we do not have to drive the original sequence with the same set of Brownian motions as only the convergence of the probability distributions of $(v^{n_j},h^{n_j})$ will matter in the first place. The next step will be to show that, on the new probability $(\tilde\Omega,\tilde{\mathcal F}, \tilde{\mathbb{P}})$, the mirror sequence converge to solution of the truncated equation. Since the 
convergence of the mirror sequence holds only in $L^{2}\left( \tilde\Omega ;C\left( [0,T];L^{2}(\mathbb{T}
^{2})^{3}\right) \right) $, we will need to resort to the weak (in probabilistic sense) version of the equation satisfied by the mirror image of  $(v^{n_j},h^{n_j})$. Let us ignore the choice of the subsequence and 
denote the mirror sequence by $((\tilde{v}^n, \tilde{h}^n), (\tilde{W}%
^{i,n})_i)_{n=1}^{\infty}$. Note that we added the corresponding set Brownian motions for each element of the sequence in the light of the discussion from above.  Then, for any test function $\varphi \in W^{3,2}(%
\mathbb{T}^2)$, we can write 

\begin{equation*}  \label{vweak}
\begin{aligned}
\langle \tilde{v}_t^n,\varphi\rangle &= \langle \tilde{v}_0^n,
\varphi\rangle + \nu \displaystyle\int_0^t\langle \tilde{v}_s^n, \Delta
\varphi \rangle ds - \displaystyle\int_0^t f_R(\tilde{a}_s^{n-1})\langle 
\tilde{v}_s^n, \mathcal{L}_{\tilde{u}_s^{n-1}}^{\star} \varphi \rangle ds - %
\displaystyle\int_0^t \langle \tilde{u}_s^n, f\hat{z} \times \varphi \rangle
ds \\
& - \displaystyle\int_0^t \langle \tilde{p}_s^{n}, \nabla
\varphi \rangle ds + \frac{1}{2}\displaystyle\sum_{i=1}^{\infty}\displaystyle%
\int_0^t \langle \tilde{v}_s^n, (\mathcal{L}_i^{\star} + \mathcal{A}%
_i^{\star})^2\varphi \rangle ds 
 - \displaystyle\sum_{i=1}^{\infty}\displaystyle\int_0^t\langle \tilde{v}%
_s^n,(\mathcal{L}_i^{\star} + \mathcal{A}_i^{\star}) \varphi \rangle d\tilde W_s^{i,n} 
\end{aligned}
\end{equation*}
\begin{equation}  \label{hweak}
\begin{aligned}
\langle \tilde{h}_t^n, \varphi \rangle &= \langle \tilde{h}_0^n,
\varphi\rangle + \eta \displaystyle\int_0^t \langle \tilde{h}_s^n, \Delta
\varphi \rangle ds - \displaystyle\int_0^t f_R(\tilde{a}_s^{n-1}) \langle
\nabla\varphi, \tilde{h}_s^{n-1} \tilde{u}_s^{n-1} \rangle ds \\
&+ \frac{1}{2}%
\displaystyle\sum_{i=1}^{\infty}\displaystyle\int_0^t\langle \tilde{h}_s^n, (%
\mathcal{L}_i^{\star})^2 \varphi \rangle ds 
 - \displaystyle\sum_{i=1}^{\infty} \displaystyle\int_0^t \langle
\tilde{h}_s^n, \mathcal{L}_i^{\star}\varphi \rangle d\tilde W_s^{i,n}. 
\end{aligned} 
\end{equation}
The next step would be to pass to the limit in (\ref{vweak}) and (\ref{hweak}) and show that each term converges to the corresponding term in the equation satisfied by truncated system. The convergence of the stochastic integrals in (\ref{vweak}) and (\ref{hweak}) poses an additional difficulty. The reason is that, even though we know that the both the integrands and the integrators (the driving Brownian motions) converge, that does not necessarily imply that the corresponding integrals converge. 
To circumvent this hurdle we make use of By Theorem 4.2 in \cite{KurtzProtter} which states that if the integrands and the integrators converge in distribution (in the original space), then the stochastic integrals converge in distributions too. Then, via the Skorokhod representation theorem, we find a mirror probability space $(\tilde\Omega,\tilde{\mathcal F}, \tilde{\mathbb{P}})$ where, \emph{by construction}, not only   
$((\tilde{v}^n, \tilde{h}^n), (\tilde{W}^{i,n})_i)_{n=1}^{\infty}$ converge, but also the corresponding stochastic integrals. We  proceed with the construction as follows:

We choose $\left\{ \varphi _{k}\right\}_{k}$ to be a countable dense set of $\mathcal{W}^{2,2}(\mathbb{T}^{2})$. By Proposition \ref{tightnessISRSW} 
and Theorem 4.2 in \cite{KurtzProtter} we can deduce that the processes 
\begin{equation*}
\begin{aligned}\{v^{n},h^{n},& \int_{{0}}^{\cdot }\left\langle
v_{1}^{n},\left( \mathcal{L}_{i}+\mathcal{A}_{i}\right) ^{\ast }\varphi
_{k_{1}}\right\rangle _{L^{2}(\mathbb{T}^{2})}dW_{s}^{i_{1},n},
\int_{{0}}^{\cdot }\left\langle v_{2}^{n},\left(
\mathcal{L}_{i}+\mathcal{A}_{i}\right) ^{\ast }\varphi _{k_{2}}\right\rangle
_{L^{2}(\mathbb{T}^{2})}dW_{s}^{i_{2},n}, \\ & \int_{{0}}^{\cdot }\left\langle
h^{n},\mathcal{L}_{i}\varphi _{k}\right\rangle
_{L^{2}(\mathbb{T}^{2})}dW_{s}^{i_{3},n}%
\text{,~~}i_{1},i_{2},i_{3},k_{1},k_{2},k_{3}=1,...\infty ,\}_{n=1}^{\infty}
\end{aligned}
\end{equation*}%
converge in distribution (possibly by extracting a subsequence of the original sequence and re-indexing it). We apply next the Skorokhod representation theorem to this (enlarged) sequence and find a probability space 
$(\tilde\Omega,\tilde{\mathcal F}, \tilde{\mathbb{P}})$
on which we can find  processes 
\begin{equation*}
\begin{aligned}\{\tilde v^{n},\tilde h^{n},& \int_{{0}}^{\cdot }\left\langle
\tilde v_{1}^{n},\left( \mathcal{L}_{i}+\mathcal{A}_{i}\right) ^{\ast }\varphi
_{k_{1}}\right\rangle _{L^{2}(\mathbb{T}^{2})}d\tilde  W_{s}^{i_{1},n},
\int_{{0}}^{\cdot }\left\langle \tilde  v_{2}^{n},\left(
\mathcal{L}_{i}+\mathcal{A}_{i}\right) ^{\ast }\varphi _{k_{2}}\right\rangle
_{L^{2}(\mathbb{T}^{2})}d \tilde  W_{s}^{i_{2},n}, \\ & \int_{{0}}^{\cdot }\left\langle
\tilde h^{n},\mathcal{L}_{i}\varphi _{k}\right\rangle
_{L^{2}(\mathbb{T}^{2})}d\tilde  W_{s}^{i_{3},n}%
\text{,~~}i_{1},i_{2},i_{3},k_{1},k_{2},k_{3}=1,...\infty ,\}_{n=1}^{\infty}
\end{aligned}
\end{equation*}%
with the same probability distributions as the corresponding elements of the original sequence and that converge to  
\begin{equation*}
\begin{aligned}(\tilde v,\tilde h,& \int_{{0}}^{\cdot }\left\langle
\tilde v_{1},\left( \mathcal{L}_{i}+\mathcal{A}_{i}\right) ^{\ast }\varphi
_{k_{1}}\right\rangle _{L^{2}(\mathbb{T}^{2})}d\tilde  W_{s}^{i_{1}},
\int_{{0}}^{\cdot }\left\langle \tilde  v_{2},\left(
\mathcal{L}_{i}+\mathcal{A}_{i}\right) ^{\ast }\varphi _{k_{2}}\right\rangle
_{L^{2}(\mathbb{T}^{2})}d \tilde  W_{s}^{i_{2}}, \\ & \int_{{0}}^{\cdot }\left\langle
\tilde h,\mathcal{L}_{i}\varphi _{k}\right\rangle
_{L^{2}(\mathbb{T}^{2})}d\tilde  W_{s}^{i_{3}}%
\text{,~~}i_{1},i_{2},i_{3},k_{1},k_{2},k_{3}=1,...\infty )
\end{aligned}
\end{equation*}%
in the corresponding product spaces as well as pathwise (possibly by extracting a suitable subsequence).   

\noindent It follows that all the estimates established for $(v^n,h^n)$ are also true for $(\tilde{v}^n,\tilde{h}^n)$. Thus, there exist a constant $\tilde{\mathcal{B}}_{3}(T,R)$ 
such that 
\begin{equation}
\tilde{\mathbb{E}}\left[ \left\vert \left\vert (\tilde v^{n},\tilde h^{n})\right\vert \right\vert
_{T,2,2}^{2}\right] \leq \tilde{\mathcal{B}}_{3}(T,R),  \label{unifestim1}
\end{equation}
which ensures that the corresponding time integrals of the terms involved
are uniformly bounded in $L^2(\tilde{\mathbb{P}})$, and, by Fatou's lemma, also  that 
\begin{equation}
\tilde{\mathbb{E}}\left[ \left\vert \left\vert (\tilde v,\tilde h)\right\vert \right\vert
_{T,2,2}^{2}\right] \leq \tilde{\mathcal{B}}_{3}(T,R),  \label{unifestim2}
\end{equation}
Similarly, we also have that  $%
(\tilde v^{n},\tilde h^{n})\in L^{p}\left( \tilde \Omega ;\mathcal{W}^{\alpha,p}\left( [0,T],
L^{2}(\mathbb{T}^{2})^{3}\right) \right) $ with $p \in(2,\infty),\alpha \in[0,\frac{1}{2})$ such
that $p\alpha \geq 1$ and there exists a constant $\mathcal{B}_{4}(p,\alpha
,T,R)$ independent of $n$ such that
\begin{equation}
\mathbb{\tilde E}\left[ \left\vert \left\vert (\tilde v^{n},\tilde h^{n})\right\vert \right\vert
_{\mathcal{W}^{\alpha, p}\left( [0,T],L^{2}(\mathbb{T}^{2})^{3}\right) }%
\right] \leq \mathcal{B}_{4}(p,T, R).  
\end{equation}
with the same control applying to the limit process 
$%
(\tilde v,\tilde h)\in L^{p}\left( \tilde \Omega ;\mathcal{W}^{\alpha,p}\left( [0,T],
L^{2}(\mathbb{T}^{2})^{3}\right) \right)$. We pass to the limit in all the terms in \eqref{vweak} and \eqref{hweak}. The stochastic terms converge by construction, therefore we only need to concentrate on the deterministic terms. Of these, the convergence of the linear terms is straightforward and relies on the convergence of $(\tilde v^{n},\tilde h^{n})$  in $L^{2}\left( \tilde\Omega ;C(\left( [0,T];L^{2}(\mathbb{T}^{2})^{3}\right) \right) $. We detail next the convergence of the nonlinear terms. For the velocity equation we show that
\begin{equation*} 
 \displaystyle\int_0^t \langle f_R(a_s^{n-1})\mathcal{L}_{\tilde{u}_s^{n-1}}\tilde{v}%
_s^{n-1} - f_R(a_s^R)\mathcal{L}_{\tilde{u}_s^{R}}\tilde{v}_s^R , \varphi \rangle ds %
\xrightarrow[n\rightarrow\infty]{} 0 \ \  \hbox{in} \
L^2(\tilde{\mathbb{P}}).  
\end{equation*}
One can split this difference as follows 
\begin{equation*}
    \begin{aligned}
    |\langle f_R(a_s^{n-1})\mathcal{L}_{\tilde{u}_s^{n-1}}\tilde{v}_s^{n-1} - f_R(a_s^R)\mathcal{L}_{\tilde{u%
}_s^{R}}\tilde{v}_s^R , \varphi \rangle| &\leq f_R(a_s^R)|\langle(\tilde{u}_s^{n-1} - 
\tilde{u}_s^R) \cdot \nabla \tilde{v}_s^{n-1}, \varphi \rangle| \\
& + f_R(a_s^R)|\langle \tilde{u}
_s^R \cdot \nabla(\tilde{v}_s^{n-1}-\tilde{v}_s^R), \varphi \rangle| \\
& + |f_R(a_s^{n-1})-f_R(a_s^R)||\langle \tilde{u}_s^{n-1} \cdot \nabla \tilde{v}_s^{n-1}, \varphi\rangle|.
    \end{aligned}
\end{equation*}
For the first term we have
\begin{equation*}
\begin{aligned}
    \mathbb{E}\left[ \displaystyle\int_0^t f_R(a_s^R)|\langle(\tilde{u}_s^{n-1} - 
\tilde{u}_s^R) \cdot \nabla \tilde{v}_s^{n-1}, \varphi \rangle|ds \right]
& \leq C(\|\varphi\|_{\infty})\mathbb{E}\left[ \displaystyle\sup_{s\in[0,t]}\|\tilde{u}_s^{n-1} - 
\tilde{u}_s^R\|_2\displaystyle\int_0^t\| \nabla\tilde{v}_s^{n-1}\|_2 ds\right] \\
& \leq C(\|\varphi\|_{2,2})\sqrt{\mathbb{E}\left[ \displaystyle\sup_{s\in[0,t]}\|\tilde{u}_s^{n-1} - 
\tilde{u}_s^R\|_2^2\right] \mathbb{E}\left[\displaystyle\int_0^t\| \nabla\tilde{v}_s^{n-1}\|_2^2 ds \right]} \\
& \leq C(t, \|\varphi\|_{2,2})\sqrt{\mathbb{E}\left[ \displaystyle\sup_{s\in[0,t]}\|\tilde{u}_s^{n-1} - 
\tilde{u}_s^R\|_2^2\right] \mathbb{E}\left[\displaystyle\sup_{s\in[0,t]}\| \tilde{v}_s^{n-1}\|_{1,2}^2 \right]} \\
& \leq C(t, \|\varphi\|_{2,2})\mathbb{E}\left[ \displaystyle\sup_{s\in[0,t]}\|\tilde{u}_s^{n-1} - 
\tilde{u}_s^R\|_2^2\right]^{1/2}
\end{aligned}
\end{equation*}
and the term on the right hand side converges to 0 in  $L^2(\tilde{\mathbb{P}})$  and all other terms are controlled uniformly in $n$. 
For the second term,
\begin{equation*}
\begin{aligned}
& \mathbb{E}\left[ \displaystyle\int_0^t f_R(a_s^R)|\langle \tilde{u}
_s^R \cdot \nabla(\tilde{v}_s^{n-1}-\tilde{v}_s^R), \varphi \rangle|ds\right]  
= \mathbb{E}\left[ \displaystyle\int_0^t f_R(a_s^R)|\langle \tilde{u}
_s^R \cdot (\tilde{v}_s^{n-1}-\tilde{v}_s^R), \nabla\varphi \rangle|ds\right] \\
& + \mathbb{E}\left[ \displaystyle\int_0^t f_R(a_s^R)|\langle (\nabla \cdot \tilde{u}
_s^R )\cdot (\tilde{v}_s^{n-1}-\tilde{v}_s^R), \varphi \rangle|ds\right] \\
& \leq \mathbb{E}\left[ \displaystyle\sup_{s\in[0,t]}\|\tilde{v}_s^{n-1}-\tilde{v}_s^R\|_2\displaystyle\int_0^tf_R(a_s^R)\|\tilde{u}_s^R \cdot \nabla\varphi\|_2ds \right] +C(\|\varphi\|_{\infty}) \mathbb{E}\left[ \displaystyle\sup_{s\in[0,t]}\|\tilde{v}_s^{n-1}-\tilde{v}_s^R\|_2 \displaystyle\int_0^tf_R(a_s^R)\|\nabla \cdot \tilde{u}_s^R\|_2ds\right] \\
& \leq C \sqrt{\mathbb{E}\left[ \displaystyle\sup_{s\in[0,t]}\|\tilde{v}_s^{n-1}-\tilde{v}_s^R\|_2^2\right]\mathbb{E}\left[ \displaystyle\int_0^t f_R(a_s^R)^2\|\tilde{u}_s^R\|_{1,2}^2\|\nabla\varphi\|_{1,2}^2\right]} \\
& +C(\|\varphi\|_{2,2}) \sqrt{\mathbb{E}\left[ \displaystyle\sup_{s\in[0,t]}\|\tilde{v}_s^{n-1}-\tilde{v}_s^R\|_2^2\right] \mathbb{E}\left[ \displaystyle\int_0^tf_R(a_s^R)^2\|\tilde{u}_s^R\|_{1,2}^2ds\right]} \\
& \leq C(t,\|\varphi\|_{2,2}, R)\mathbb{E}\left[ \displaystyle\sup_{s\in[0,t]}\|\tilde{v}_s^{n-1}-\tilde{v}_s^R\|_2^2\right]^{1/2} \ \ \xrightarrow[n\rightarrow \infty]{}0. 
\end{aligned}
\end{equation*}
For the third term,
\begin{equation*}
\begin{aligned}
   & \mathbb{E}\left[\displaystyle\int_0^t|f_R(a_s^{n-1})-f_R(a_s^R)||\langle \tilde{u}_s^{n-1} \cdot \nabla \tilde{v}_s^{n-1}, \varphi\rangle|ds\right]\\  
   \leq & C(\|\varphi\|_{\infty}) \sqrt{\mathbb{E}\left[\displaystyle\int_0^t|f_R(a_s^{n-1})-f_R(a_s^R)|^2ds\right]\mathbb{E}\left[\displaystyle\int_0^t\|\tilde{u}_s^{n-1} \cdot \nabla \tilde{v}_s^{n-1}\|_2^2ds\right]} \\
   & \leq C(t,\|\varphi\|_{2,2})\sqrt{\mathbb{E}\left[\displaystyle\sup_{s\in[0,t]}\|a_s^{n-1}-a_s^R\|_2\displaystyle\int_0^t\|a_s^{n-1}-a_s^R\|_{2,2}ds\right]\mathbb{E}\left[\displaystyle\sup_{s\in[0,t]}\|a_s^{n-1}\|_{1,2}^2\right]} \\
   & \leq C(t,\|\varphi\|_{2,2})\mathbb{E}\left[\displaystyle\sup_{s\in[0,t]}\|a_s^{n-1}-a_s^R\|_2^2\right]^{1/4}\mathbb{E}\left[ \displaystyle\int_0^t(\|a_s^{n-1}\|_{2,2}^2 + \|a_s^R\|_{2,2}^2) ds\right]^{1/4} \\
   & \leq \tilde{C}(t,\|\varphi\|_{2,2})\mathbb{E}\left[\displaystyle\sup_{s\in[0,t]}\|a_s^{n-1}-a_s^R\|_2^2\right]^{1/4} \ \ \xrightarrow[n\rightarrow \infty]{}0.
\end{aligned}
\end{equation*}
\noindent Note that $\mathbb{E}\left[ \displaystyle\int_0^t\|a_s^R\|_{2,2}^2 ds\right] < \infty$ by a direct application of the Fatou lemma. 
With similar arguments, the nonlinear term in the height equation $\eqref{hweak}$ 
converges as requested:
\begin{equation*}\displaystyle\int_0^t \langle f_R(a_s^{n-1}) \nabla \cdot (\tilde{h}_s^{n-1}\tilde{u}%
_s^{n-1}) - f_R(a_s^R)\nabla \cdot(\tilde{h}_s^{R}\tilde{u}_s^R), \varphi \rangle \rangle ds %
\xrightarrow[n\rightarrow\infty]{} 0 \ \ \ \hbox{in} \ \ L^2(\tilde{\mathbb{P}}).
\end{equation*}

\noindent We have constructed a weak (in PDE sense) solution of the SRSW system, as we have chosen the set of test functions  $(\varphi _{k})_{k}$ to be a countable dense set of $\mathcal{W}^{2,2}(\mathbb{T}^{2})$. Since $(\tilde v, \tilde h)$ has the right amount of smoothness, this weak solution is also strong (in PDE sense). However, 
$(\tilde v, \tilde h)$ is constructed on a different probability space than the original one. We apply next the Yamada-Watanabe theorem (see, e.g. Theorem 2.1 in \cite{Roeckner}) 
to justify that the existence of the solution on this different probability space together with the pathwise unique of the truncated equation implies the existence of a (unique) solution of the truncated equation on the original space. We have constructed a weakly continuous solution $a^R \in L^2\left( \Omega, L^{\infty}\left([0,T]; \mathcal{W}^{1,2}(\mathbb{T}^2)^3\right)\right)$. From Lemma \ref{lemma:globsol1} we can deduce that $\mathbb{E}\left [\left(\|a_t^R\|_{1,2}^2 - \|a_s^R\|_{1,2}^2\right)^4\right] \leq C(t-s)^2$, and therefore by the Kolmogorov-Čentsov criterion, the map $t \rightarrow \|a_t^R\|_{1,2}^2$ is continuous. Hence $a^R \in L^2\left( \Omega, C\left(0,T; \mathcal{W}^{1,2}(\mathbb{T}^2)^3\right)\right)$. 

The proof of the claim is now complete. 

\section{Global solution with positive probability}
Let $(a, \mathcal{T})$ be a maximal solution of the SRSW system and recall that $    \tau_R = \displaystyle\inf_{t \geq 0}\{\|a_t\|_{1,2} >R\} $. The following technical lemma gives the main estimate for showing the global solution property.    
\begin{lemma}\label{lemma:globsol1}
Let $(a, \mathcal{T})$ be a maximal solution of the SRSW system. Then there exist some positive constants $C_i, i=1,3$, independent of $R$ such that 
\begin{equation*}
    \|a_{t\wedge\tau_R}\|_{1,2}^2 = \|a_0\|_{1,2}^2 + \displaystyle\int_0^{t\wedge\tau_R}\tilde{F}(a_s)ds + \displaystyle\sum_{i=1}^{\infty}\displaystyle\int_0^{t\wedge\tau_R}\tilde{G}_i(a_s)dW_s^i
\end{equation*}
where $\tilde{F}(a_s)$ and $\tilde{G}_i(a_s)$ are processes such that
\begin{equation*}
\begin{aligned}
    & |\tilde{F}(a_s)|^2 \leq C_1\|a_s\|_{1,2}^6 - C_2\|a_s\|_{1,2}^2 \\
   &  \displaystyle\sum_{i=1}^{\infty}|\tilde{G}_i(a_s)|^2 \leq C_2\|a_s\|_{1,2}^2.
\end{aligned}
\end{equation*}
\end{lemma}

\noindent The proof of this lemma is provided in the Appendix. 

\begin{proposition}
Let $(a, \mathcal{T})$ be a maximal solution of the SRSW system. Then $\tau_R>0$ $\mathbb{P}$-a.s. for any  $R> \|a_0\|_{1,2}$. In particular $\mathcal{T}>0$ $\mathbb{P}$-a.s.
\end{proposition}
\noindent\textbf{Proof.} From Lemma \ref{lemma:globsol1} and the Burkholder-Davis-Gundy inequality we deduce that
\begin{equation*}
    \mathbb{E}\left[ |\|a_{t\wedge\tau_R}\|_{1,2}^2 - \|a_0\|_{1,2}^2|\right]\leq tR^6 + \sqrt{t}R^2.
\end{equation*}
Note that on the set $\{ \tau_R < t\}$ we have $\|a_{t\wedge\tau_R}\|_{1,2} =R$. It follows that
\begin{equation*}
    \begin{aligned}
    (R^2 - \|a_0\|_{1,2}^2)\mathbb{P}(\tau_R < t)& = \mathbb{E}\left[ |\|a_{t\wedge\tau_R}\|_{1,2}^2 - \|a_0\|_{1,2}^2|\mathds{1}_{\{\tau_R < t\}}\right] \\
    & \leq \mathbb{E}\left[ |\|a_{t\wedge\tau_R}\|_{1,2}^2 - \|a_0\|_{1,2}^2|\right] \\
    & \leq tR^6+\sqrt{t}R^2. 
    \end{aligned}
\end{equation*}
Hence
\begin{equation*}
    \mathbb{P}(\tau_R<t)\leq\frac{tR^6+\sqrt{t}R^2}{R^2-\|a_0\|_{1,2}^2}.
\end{equation*}
Then
\begin{equation*}
    \displaystyle\lim_{t\rightarrow 0}\mathbb{P}(\tau_R < t) = 0
\end{equation*}
and
\begin{equation*}
    \mathbb{P}(\tau_R = 0) = \displaystyle\bigcap_{n>0}\mathbb{P}\left(\tau_R< \frac{1}{n}\right) = \displaystyle\lim_{n\rightarrow\infty}\mathbb{P}\left(\tau_R < \frac{1}{n}\right) = 0.
\end{equation*}
Hence $\tau_R > 0$, $\mathbb{P}$-a.s. and therefore also $\mathcal{T} \geq \tau_R$ is strictly positive $\mathbb{P}$ almost surely. 

\begin{proposition}
Let $(a, \mathcal{T})$ be a maximal solution. Then there exists a positive
constant $C$ such that, if $\|a_0\|_{1,2} < C$ then $\mathbb{P}(\mathcal{T}=\infty)>0$. In
other words, if the initial condition is sufficiently small, then the equation
has a global solution.
\end{proposition}

\noindent\textbf{Proof} Using the notation in Lemma \ref{lemma:globsol1}, define
\begin{equation*}
A(a_s)= 
\begin{cases}
\frac{\tilde{F}(a_s)}{\|a_s\|_{1,2}^2}, & \text{if \ \ $a_s \neq 0$}. \\ 
0, & \text{if \ \ $a_s=0$}.%
\end{cases}
\end{equation*}
\begin{equation*}
 B_i(a_s)=  
\begin{cases}
\frac{\tilde{G}_i(a_s)}{\|a_s\|_{1,2}^2}, & \text{if \ \ $a_s \neq 0$}. \\ 
0, & \text{if \ \ $a_s=0$}.%
\end{cases}
\end{equation*}
We deduce from Lemma \ref{lemma:globsol1} that 
\begin{equation*}  \label{main2}
\begin{aligned} 
\|a_{t\wedge\tau_R}\|_{1,2}^2 & = \|a_0\|_{1,2}^2 + \displaystyle\int_0^{t\wedge \tau_R}
A(a_s)\|a_s\|_{1,2}^2 ds +
\displaystyle\sum_{i=1}^{\infty}\displaystyle\int_0^{t\wedge \tau_R}
B_i(a_s)\|a_s\|_{1,2}^2 dW_s^i. 
\end{aligned}
\end{equation*}
This implies that 
\begin{equation*}
\begin{aligned} 
\|a_{t\wedge\tau_R}\|_{1,2}^2 = \|a_0\|_{1,2}^2\exp \left(
\displaystyle\int_0^{t\wedge \tau_R} A(a_s)ds + M_{t\wedge \tau_R} -\frac{1}{2}[M]_{t\wedge \tau_R}\right) \end{aligned}
\end{equation*}
where $M$ is the local martingale defined (for $t\geq 0$) as 
\begin{equation*}
M_t = \displaystyle\sum_{i=1}^{\infty}\displaystyle\int_0^tB_i(a_s)dW_s^i
\end{equation*}
with quadratic variation given by
\begin{equation*}
[M]_t = \displaystyle\sum_{i=1}^{\infty}\displaystyle\int_0^tB_i(a_s)^2ds.
\end{equation*}
Moreover, since 
\begin{equation*}
    \displaystyle\sum_{i=1}^{\infty}|\tilde{G}_i(a_s)|^2 \leq C_3\|a_s\|_{1,2}^2
\end{equation*}
we have that
\begin{equation*}
    \displaystyle\sum_{i=1}^{\infty}|B_i(a_s)|^2 \leq C_3.
\end{equation*}
It follows that $M$ is a square integrable martingale with quadratic
variation $[M]_t \leq C_3$.
In particular, by Novikov condition, 
$\exp \left(M_t -\frac{1}{2}[M]_t\right)$
is a martingale and therefore 
$\mathbb{E}\left[ \exp\left(M_{t\wedge\tau_R} -\frac{1}{2}[M]_{t\wedge\tau_R}\right) \right] = 1.$
Next we have from Lemma \ref{abstractlemmaexistence} that 
\begin{equation*}
\tilde{F}(a_s) \leq c_1 \|a_s\|_{1,2}^6 - c_2 \|a_s\|_{1,2}^2
\end{equation*}
hence 
\begin{equation*}
A(a_s) \leq c_1\|a_s\|_{1,2}^4 - c_2.
\end{equation*}
Choose 
\begin{equation*}
\|a_0\|_{1,2}^2 < \left(\frac{c_2}{c_1}\right)^{1/4} =: C
\end{equation*}
and define 
\begin{equation*}
\tau_C :=\displaystyle\inf_{t} \{ \|a_t\|_{1,2}^2 \geq C \}.
\end{equation*}
Then 
\begin{equation*}
\begin{aligned} \mathbb{E}\left[ \|a_{t\wedge \tau_C}\|_{1,2}^2\right] & =
\|a_0\|_{1,2}^2\mathbb{E}\left[\exp \left(
\displaystyle\int_0^{t\wedge\tau_C}A(a_s)ds + M_{t\wedge\tau_C} -
\frac{1}{2}[M]_{t\wedge\tau_C}\right) \right] \\ & <
\|a_0\|_{1,2}^2\mathbb{E}\left[\exp \left( M_{t\wedge\tau_C} -
\frac{1}{2}[M]_{t\wedge\tau_C}\right) \right] \\ & < C. \end{aligned}
\end{equation*}
Now 
\begin{equation*}
\begin{aligned} \mathbb{P}\left( \tau_C < \infty \right) & =
\displaystyle\bigcap_N\mathbb{P}(\tau_C > N) \\ & =
\displaystyle\lim_{N\rightarrow\infty}\mathbb{P}(\|a_{N\wedge\tau_C}%
\|_{1,2}^2 <C) \\ 
& =\displaystyle\lim_{N\rightarrow\infty}\left(1-\mathbb{P}\left(\|a_{N\wedge\tau_C}%
\|_{1,2}^2 \geq C\right)\right) \\ 
& \leq 1 - \frac{\|a_0\|_{1,2}^2}{C} \end{aligned}
\end{equation*}
since we have 
\begin{equation*}
\mathbb{P}(\|a_{N\wedge\tau_C}\|_{1,2}^2 \geq C) \leq \frac{\mathbb{E}%
[\|a_{N\wedge\tau_C}\|_{1,2}^2]}{C} \leq \frac{\|a_0\|_{1,2}^2}{C} < 1.
\end{equation*}
It follows that $\mathbb{P}(\tau_R =\infty) > 0$ hence the claim. 

\section{Analytical properties of the approximating system}
\label{section:analyticalapprox}

\subsection{Relative compactness}
\label{subsection:relcompactness} 
We define the following processes, to
shorten the notation: 
\begin{equation*}
\begin{aligned} & X_t^{v^n} := v_0^n + \displaystyle\int_0^t\left(\nu \Delta
v_s^n + P_s^{n-1,n}(v_s^n)\right)ds \\ & Y_t^{v^n} := \displaystyle\int_0^t
\displaystyle\sum_{i=1}^{\infty}[(\mathcal{L}_i +
\mathcal{A}_i)v_s^n]dW_s^{i,n}\\ & X_t^{h^n} := h_0^n +
\displaystyle\int_0^t\left(\eta \Delta h_s^n + Q_s^{n-1,n}(h_s^n)\right)ds
\\ & Y_t^{h^n} := \displaystyle\int_0^t
\displaystyle\sum_{i=1}^{\infty}[\nabla \cdot (\xi_ih_s^n)]dW_s^{i,n}.
\end{aligned}
\end{equation*}

\noindent \textbf{\textit{Proof of Theorem \ref{psnsSRSW}}} The existence and uniqueness of the
solution of the system follows directly from Theorem \ref{theoremrozovskii}.
The
control (\ref{c1}) holds true from the same theorem and the fact that all
coefficients are the same with the exception of the forcing
terms, which are bounded uniformly in $n$, as we show below. Let 
\begin{equation*}
F_s^{n-1} = F_s^{n-1,u} + F_s^{n-1,h}:=f_{R}(a_{s}^{n-1}) \left( u_{s}^{n-1} \cdot \nabla v_{s}^{n-1} + \nabla \cdot (
h_{s}^{n-1}u_{s}^{n-1})\right).
\end{equation*}%
The $L^2$ norm of the first term can be controlled using the truncation and Ladyzhenskaya's inequality, as follows\footnote{Note that $C$ can be different at each line.}
\begin{equation*}
\begin{aligned}
\displaystyle\int_0^t\|F_s^{n-1,u} \|_2^2 ds& =\displaystyle\int_0^t f_{R}(a_{s}^{n-1})  \|u_{s}^{n-1} \cdot \nabla v_{s}^{n-1}\|_2^2ds  \leq C\displaystyle\int_0^tf_R(a_s^{n-1})\|u_s^{n-1}\|_4^2\|\nabla v_s^{n-1}\|_4^2ds \\
& \leq CR^3\displaystyle\int_0^t\|v_s^{n-1}\|_{2,2}ds  \leq CR^3\sqrt{t}\sqrt{\displaystyle\int_0^t\|v_s^{n-1}\|_{2,2}^2ds} \\
& \leq C\sqrt{\tilde{C}}R^3\sqrt{t} \leq C_1(R,t).
\end{aligned}
\end{equation*}
Similarly, using Lemma \ref{abstractlemmaexistence} from Appendix we have  that 
\begin{equation*}
    \begin{aligned}
   \displaystyle\int_0^t\|F_s^{n-1,h} \|_2^2 ds &:= \displaystyle\int_0^tf_R(a_s)\|\nabla \cdot (h_s^{n-1}u_s^{n-1}))\|_2^2 \leq CR^3\displaystyle\int_0^t(\|h_s^{n-1}\|_{2,2} + \|u_s^{n-1}\|_{2,2})ds \\
   & \leq 2C\sqrt{\tilde{C}}R^3\sqrt{t} \leq C_2(R,T).
    \end{aligned}
\end{equation*}
Summing up and using an inductive argument we deduce that there exists a constant $C$ which is independent of $n$ such that
\begin{equation*}
    \begin{aligned}
  \mathbb{E}\left[ \|a_s^n\|_{t,2,2}^2\right] &\leq N\left( \|a_0\|_{1,2}^2 + \mathbb{E}\left[ \displaystyle\int_0^t \|F_s^{n-1}\|_2^2ds\right] \right)  
         \leq N \left( \|a_0\|_{1,2}^2 + CR^3\sqrt{t}\right) \leq C(R,t).
    \end{aligned}
\end{equation*}
For an arbitrary $p>2$, we can deduce that there exists a constant $\widetilde{\mathcal{D}}_{p}(T,R)$  such that 
\[
\mathbb E [\|a^n\|_{T,2,2}^p ]\le N \left(\|a_0\|_{1,2}^p + \widetilde{\mathcal{D}}_{p}(T,R)
\sqrt{\mathbb E [ \|a^{n-1}\|_{T,2,2}^p]} \right) 
\]
The result follows with an argument similar to the one used above. 
For the second part, recall that
\begin{equation*}
\|a_t^n\|_{\mathcal{W}^{\beta, p}(0, T ; L^2(\mathbb{T}^2))}^p:=\int_{0}^{T}\left\|a_{t}^n\right\|_{L^2(\mathbb{T}^2)}^{p} d t+\int_{0}^{T} \int_{0}^{T} \frac{\left\|a_{t}^n-a_{s}^n\right\|_{L^2(\mathbb{T}^2)}^{p}}{|t-s|^{1+\beta p}} d t d s.
\end{equation*} 
We show that there exists a constant $C=C\left(
T,R\right) $ independent of $n$ such that
\begin{equation*}
E\left[ \left\|a_{t}^{n}-a_{s}^{n}\right\| _{L^2(\mathbb{T}^2)}^{p}%
\right] \leq C |t-s|^{p/2}.
\end{equation*}
We have
\begin{equation*}
    \begin{aligned}
    X_t^{v^n}-X_s^{v^n} =    \displaystyle\int_s^t P_r^{n-1,n}(v_r^n) dr + \displaystyle\int_s^t\nu\Delta v_r^ndr.
    \end{aligned}
\end{equation*}

Then
\begin{equation*}
    \begin{aligned}
    \mathbb{E}\left[\|X_t^{v^n}-X_s^{v^n}\|_2^p \right]& \leq \mathbb{E}\left[\left( \displaystyle\int_s^t \|P_r^{n-1,n}(v_r^n)\|_2dr\right)^p + \left( \displaystyle\int_s^t\|\Delta v_r^n\|_2 dr\right)^p \right]\\
    &\leq  (t-s)^{p}\mathbb{E}\left[\displaystyle\sup_{r\in[s,t]} \|P_r^{n-1,n}(v_r^n)\|_2^p \right] + \mathbb{E}\left[(t-s)^{p/2}\left(\displaystyle\int_0^T\|\nu\Delta v_r^n\|_2^2dr\right)^{p/2}\right] \\
    & \leq (t-s)^{p}\mathbb{E}\left[\displaystyle\sup_{r\in[0,T]} \|P_r^{n-1,n}(v_r^n)\|_2^p \right] + (t-s)^{p/2}\mathbb{E}\left[ \|v_r^n\|_{T,2,2}^p\right] \\
    & \leq C(t-s)^{p/2}.
    \end{aligned}
\end{equation*}
For the stochastic terms we apply the Burkholder-Davis-Gundy inequality to obtain
\begin{equation*}
    \begin{aligned}
    \mathbb{E}\left[\|Y_t^{v^n}-Y_s^{v^n}\|_2^p \right] &\leq \mathbb{E}\left[\left|\displaystyle\int_s^t \displaystyle\sum_{i=1}^{\infty}\langle v_r^n, (\mathcal{L}_i + \mathcal{A}_i)v_r^n\rangle dW_r^{i,n}\right|^{p}\right] \\
& \leq C(p) \mathbb{E}\left[\displaystyle\int_s^t
\displaystyle\sum_{i=1}^{\infty}|\langle v_r^n, (\mathcal{L}_i + \mathcal{A}_i)v_r^n\rangle|^2 dr\right]^{p/2} \\
& \leq C(p) \mathbb{E}\left[ \displaystyle\int_s^t\|
v_r^n\|_{2}^2 \displaystyle\sum_{i=1}^{\infty}\|(\mathcal{L}_i + \mathcal{A}_i)v_r^n\|_{2}^2dr\right]^{p/2} \\
& \leq C(p) \mathbb{E}\left[ \displaystyle\int_s^t\|
v_r^n\|_{2}^2 \|v_r^n\|_{1,2}^2dr\right]^{p/2} \\
& \leq C(p)(t-s)^{p/2}\mathbb{E}\left[ \displaystyle\sup_{r\in[s,t]}\|v_r^n\|_{1,2}^{2p}\right] \\
& \leq C(p,T)(t-s)^{p/2}. 
\end{aligned}
\end{equation*}
With similar arguments
\begin{equation*}
    \begin{aligned}
    \mathbb{E}\left[\|X_t^{h^n}-X_s^{h^n}\|_2^p \right] & \leq C(t-s)^{p/2}. 
    \end{aligned}
\end{equation*}
and
\begin{equation*}
    \begin{aligned}
    \mathbb{E}\left[\|Y_t^{h^n}-Y_s^{h^n}\|_2^p \right] & \leq C(t-s)^{p/2}. 
    \end{aligned}
\end{equation*}

\begin{proposition}
The approximating sequence is relatively compact in the space 
\begin{equation*}
C\left( [0,T], L^2(\mathbb{T}^2)^3\right).
\end{equation*}
\end{proposition}
\noindent\textbf{Proof.} 
By a standard Arzela-Ascoli argument (see e.g. \cite{SIM}), the following compact embedding holds
\begin{equation*}
L^{\infty}\left( [0,T] , \mathcal{W}^{1,2}(\mathbb{T}^2)^3\right) \cap \mathcal{W}^{\beta,p}\left(0,T; L^2(\mathbb{T}^2)^3\right) \hookrightarrow C\left([0,T], L^2(\mathbb{T}^2) \right). 
\end{equation*}
This implies that the intersection $B^N:=B^1(0,N)\cap B^2(0,N)$ of any two balls $B^1(0,N) \in L^{\infty}\left( [0,T] , \mathcal{W}^{1,2}(\mathbb{T}^2)^3\right)$ and $B^2(0,N) \in \mathcal{W}^{\beta,p}\left(0,T; L^2(\mathbb{T}^2)^3\right)$ is a compact set in the space $C\left([0,T], L^2(\mathbb{T}^2) \right)$. Observe that
\begin{equation*}
    \mathbb{P}\left( a^n \notin B^1(0,N)\right) = \mathbb{P}\left( \displaystyle\sup_{s\in[0,T]}\|a_s^n\|_{1,2}>N\right) \leq \frac{\mathbb{E}\left[\displaystyle\sup_{s\in[0,T]}\|a_s^n\|_{1,2}^2 \right]}{N^2}
\end{equation*}
\begin{equation*}
    \mathbb{P}\left( a^n \notin B^2(0,N)\right) = \mathbb{P}\left( \|a_s^n\|_{\mathcal{W}^{\beta,p}}>N\right) \leq \frac{\mathbb{E}\left[ \|a_s^n\|_{\mathcal{W}^{\beta,p}}^p\right]}{N^p}.
\end{equation*}
Hence
\begin{equation*}
\begin{aligned}
    \displaystyle\lim_{N\rightarrow\infty}\displaystyle\sup_n\mathbb{P}\left( a^n \notin B^N\right) \leq \displaystyle\lim_{N\rightarrow\infty} \frac{\displaystyle\sup_n\mathbb{E}\left[\displaystyle\sup_{s\in[0,T]}\|a_s^n\|_{1,2}^2\right]}{N^2} + \frac{\displaystyle\sup_n\mathbb{E}\left[ \|a_s^n\|_{\mathcal{W}^{\beta,p}}^p\right]}{N^p} = 0. 
\end{aligned}
\end{equation*}
This justifies the relative compactness of the distribution of $a^n$, that is the tightness of the process $a^n$, provided $\displaystyle\sup_n\mathbb{E}\left[\displaystyle\sup_{s\in[0,T]}\|a_s^n\|_{1,2}^2\right] < \infty$ and $\displaystyle\sup_n\mathbb{E}\left[ \|a_s^n\|_{\mathcal{W}^{\beta,p}}^p\right] < \infty$. These last two statements are true due to Theorem \ref{psnsSRSW} which was proven above. 
\section{Appendix}
\begin{lemma}\label{abstractlemmauniqueness}
Let $(X^1,\tau^1), (X^2,\tau^2)$ be two local solutions of the SRSW system, and $$\bar{X}:=X^1-X^2, \ \ \ \tau^{1,2} = \tau^1\wedge\tau^2, \ \ \ a^i:=(X^i, Y^i), \ \ \ \bar{a}:= a^1 - a^2,$$ 
$$Q({\bar{Y}},\bar{X}):=f_R(a^1)\nabla \cdot (Y^1 X^1) - f_R(a^2)\nabla \cdot(Y^2X^2)$$ where $\bar{Y}:=Y^1-Y^2$ and $Y$ depends linearly on $X$. 
Then there exists $\zeta>0$ and $C(\zeta,R)$ such that for $|\alpha|\leq k$
\begin{equation*}
\begin{aligned}
	|\langle \partial^{\alpha}\bar{a}, \partial^{\alpha}Q({\bar{Y}},\bar{X})\rangle| &\leq \zeta\|\bar{a}\|_{k+1,2}^2 + C(\zeta,R)\|Z\|\|\bar{a}\|_{k,2}^2
\end{aligned}
\end{equation*}
with 
\begin{equation*}
   \|Z\|:= C(\|a^1\|_{k,2}^4 + \|a^2\|_{k,2}^4)
\end{equation*}
\end{lemma}
\noindent\textbf{Proof}
We use the decomposition
\begin{equation*}
\begin{aligned}
Q(\bar{Y},\bar{X}) &= f_R(a^1)\nabla \cdot (Y^1 \bar{X}) + f_R(a^2)\nabla\cdot(\bar{Y} X^2) + |f_R(a^1) - f_R(a^2)| \nabla\cdot(Y^1X^2) \\
& := T_1 + T_2 + T_3.
\end{aligned}
\end{equation*}
\begin{itemize}
\item We have $\nabla \cdot (XY) = X \cdot \nabla Y + Y(\nabla \cdot X) = \mathcal{L}_X Y + \mathcal{D}_X Y$ for any vector $X$ and scalar $Y$. 
\item For $\mathcal{L}_X Y$ use the fact that
\begin{equation*}
\begin{aligned}
|\langle \partial^{\alpha}\bar{a}, \partial^{\alpha}(X \cdot \nabla Y)\rangle| & = |\langle \partial^{\alpha + 1}\bar{a}, \partial^{\alpha - 1}(X \cdot \nabla Y)\rangle|	\\
& \leq C \|\partial^{\alpha + 1} \bar{a}\|_2\|\partial^{\alpha-1}(X \cdot \nabla Y)\|_{2} \\
& \leq C \|\partial^{\alpha + 1} \bar{a}\|_2 \displaystyle\sum_{\beta \leq \alpha -1}C\|\partial^{\beta} X\|_4\|\partial^{\alpha-\beta} Y\|_4 \\
& \leq C \|\partial^{\alpha + 1} \bar{a}\|_2\|\partial^{\beta} X\|_2^{1/2}\|\partial^{\beta +1}X\|_2^{1/2}\|\partial^{\alpha-\beta} Y\|_2^{1/2}\|\partial^{\alpha - \beta +1} Y\|_2^{1/2}\\
& \leq C \|\bar{a}\|_{k+1,2}\|X\|_{k-1,2}^{1/2}\|X\|_{k,2}^{1/2}\|Y\|_{k,2}^{1/2}\|Y\|_{k,2}^{1/2} \\
& \leq C \|\bar{a}\|_{k+1,2}\|X\|_{k,2}\|Y\|_{k,2} \\
& \leq \frac{\zeta}{2} \|\bar{a}\|_{k+1,2}^2 + \frac{C_1(\zeta,R)}{2} \|X\|_{k,2}^2\|Y\|_{k,2}^2. 
\end{aligned}	
\end{equation*}

\item Similarly, for $\mathcal{D}_X Y$ use the fact that 
\begin{equation*}
    \begin{aligned}
    |\langle \partial^{\alpha}\bar{a}, \partial^{\alpha}(Y (\nabla \cdot X))\rangle| & = |\langle \partial^{\alpha + 1} \bar{a}, \partial^{\alpha - 1}(Y (\nabla \cdot X))\rangle| \\
    & \leq C \|\partial^{\alpha +1} \bar{a}\|_2 \displaystyle\sum_{\beta\leq\alpha - 1} C \|\partial^{\beta} Y\|_4\|\partial^{\alpha - \beta} X\|_4 \\
    & \leq \frac{\zeta}{2} \| \bar{a}\|_{k+1,2}^2 + \frac{C_2(\zeta,R)}{2} \|Y\|_{k,2}^2\| X\|_{k,2}^2.
    \end{aligned}
\end{equation*}
\item Summing up we have
\begin{equation*}
    \begin{aligned}
    |\langle\partial^{\alpha}\bar{a}, \partial^{\alpha}(\nabla \cdot (YX)) \rangle| \leq \zeta\|\bar{a}\|_{k+1,2}^2 + C(\zeta,R)\|Y\|_{k,2}^2\|X\|_{k,2}^2.
    \end{aligned}
\end{equation*}
\item Now apply this for $T_i, i=1,2,3:$
\begin{equation*}
\begin{aligned}
|\langle \partial^{\alpha}\bar{a}, \partial^{\alpha}T_1\rangle| & 
& \leq \frac{\zeta}{3}\|\bar{a}\|_{k+1,2}^2 + \frac{\tilde{C}(\zeta,R)}{3}\|Y^1\|_{k,2}^2\|\bar{X}\|_{k,2}^2 \\
\end{aligned}
\end{equation*}
\begin{equation*}
\begin{aligned}
|\langle \partial^{\alpha}\bar{a}, \partial^{\alpha}T_2\rangle| & 
& \leq \frac{\zeta}{3}\|\bar{a}\|_{k+1,2}^2 + \frac{C(\zeta,R)}{3}\|X^2\|_{k,2}^2\|\bar{Y}\|_{k,2}^2
\end{aligned}
\end{equation*}
\begin{equation*}
\begin{aligned}
|\langle \partial^{\alpha}\bar{a}, \partial^{\alpha}T_3\rangle| & \leq C|f_R(a_1)-f_R(a_2)|\|\bar{a}\|_{k+1,2}\|Y^1\|_{k,2}\|X^2\|_{k,2}\\
& \leq C\|\bar{a}\|_{k,2}\|\bar{a}\|_{k+1,2}\|Y^1\|_{k,2}\|X^2\|_{k,2}\\
& \leq \frac{\zeta}{3}\|\bar{a}\|_{k+1,2}^2 + \frac{\zeta}{3}\|Y^1\|_{k,2}^2\|X^2\|_{k,2}^2\|\bar{a}\|_{k,2}^2.
\end{aligned}
\end{equation*}
\end{itemize}
Define
\begin{equation*} 
\tilde{Z}:= C\left(\|Y^1\|_{k,2}^2 + \|X^2\|_{k,2}^2 + \|Y^1\|_{k,2}^2\|X^2\|_{k,2}^2\right) \leq C(\|Y^1\|_{k,2}^4 + \|X^2\|_{k,2}^4) \leq C(\|a^1\|_{k,2}^4 + \|a^2\|_{k,2}^4)
\end{equation*}
and 
\begin{equation*}
    \|Z\|:=C(\|a^1\|_{k,2}^4 + \|a^2\|_{k,2}^4).
\end{equation*}
Then 
\begin{equation*}
\begin{aligned}
	|\langle \partial^{\alpha}\bar{a}, \partial^{\alpha}Q({\bar{Y}},\bar{X})\rangle| &\leq \zeta\|\bar{a}\|_{k+1,2}^2 + C(\zeta,R)\|Z\|\|\bar{a}\|_{k,2}^2.
\end{aligned}
\end{equation*}

\begin{lemma}\label{abstractlemmaexistence}
\label{aprioriineqoriginal} The following statements are true: \newline
a. There exist some constants $C_1, C_2$ such that for any two vectors $X$ and $Y$  such that $Y=\epsilon X + \mathcal{R}$ and for any multi-index $\alpha$ such that $|\alpha| < k$, we have 
\begin{equation*}
|\langle \partial^{\alpha} Y_t, \partial^{\alpha}(\mathcal{L}_{X_t}Y_t)\rangle | \leq
C_1\|Y_t\|_{k+1,2}^2 + C_2\|Y_t\|_{k,2}^6.
\end{equation*}
b. There exist some constants $C_3, C_4, C_5$ such that for any scalar $Y$ and vector $X$ we have
\begin{equation*}
|\langle \Delta Y_t, \nabla \cdot (Y_tX_t)\rangle| \leq C_3 \|\Delta
Y_t\|_2^2 + C_4\|\Delta X_t\|_2^{2} + C_5(\|Y_t\|_{1,2}^6 +
\|X_t\|_{1,2}^6)
\end{equation*}
\end{lemma}

\noindent\textbf{\textit{Proof}} a. Using Agmon's, Hölder's, and Young's inequalities and the fact that $v = \epsilon u + \mathcal{R}$, we have: 
\begin{equation*}
\begin{aligned} |\langle \partial^{\alpha} Y_t, \partial^{\alpha}(\mathcal{L}_{X_t}Y_t)\rangle| &
= |\langle \partial^{\alpha+1} Y_t, \mathcal{L}_{X_t}Y_t\rangle| \leq C(\|Y_t\|_{k+1,2}
\|X_t\|_{\infty}\|\partial^{\alpha} Y_t\|_2) \\ 
& \leq C(\|
X_t\|_{k+1,2}^{1/2}\|X_t\|_{k,2}^{1/2}\|Y_t\|_{k+1,2}\|Y_t\|_{k,2}) \\ 
& \leq
C(\|Y_t\|_{k+1,2}^{1/2}\|Y_t\|_{k,2}^{1/2}\| Y_t\|_{k+1,2}\|Y_t\|_{k,2})\\ 
& \leq C(\| Y_t\|_{k+1,2}^{3/2}\|Y_t\|_{k,2}^{3/2})
\\ 
& \leq C_1\|Y_t\|_{k+1,2}^2 + C_2\|Y_t\|_{k,2}^6. \end{aligned}
\end{equation*}

b. Using Hölder's, Ladyzhenskaya's, and Young's  inequalities we have that 
\begin{equation*}
\begin{aligned} |\langle \Delta Y_t, \nabla \cdot (Y_tX_t)\rangle| & =
|\langle \Delta Y_t, Y_t\nabla \cdot X_t + X_t\nabla Y_t \rangle| \leq
\|\Delta Y_t\|_2\|Y_t\nabla X_t\|_2 + \|\Delta Y_t\|_2\|X_t \nabla \cdot
Y_t\|_2 \\ 
& \leq C(\|\Delta Y_t\|_2\|Y_t\|_4\|\nabla X_t\|_4 + \|\Delta
Y_t\|_2\|X_t\|_4 \|\nabla Y_t\|_4) \\ 
& \leq C(\|\Delta
Y_t\|_2\|Y_t\|_{1,2}\|\nabla X_t\|_2^{1/2}\|\Delta X_t\|_2^{1/2} + \|\Delta
Y_t\|_2 \|X_t\|_{1,2}\|\nabla Y_t\|_2^{1/2}\|\Delta Y_t\|_2^{1/2}) \\ 
& \leq
C(\|\Delta Y_t\|_2\|Y_t\|_{1,2}\|\nabla X_t\|_2^{1/2}\|\Delta X_t\|_2^{1/2} +
\|\Delta Y_t\|_2^{3/2} \|X_t\|_{1,2}\| Y_t\|_{1,2}^{1/2} )\\
& \leq C( \|\Delta Y_t\|_{2}^{4/3}\|Y_t\|_{1,2}^{4/3}\|X_t\|_{1,2}^{2/3} + \|\Delta X_t\|_2^2 + \|X_t\|_{1,2}^6 + \|\Delta Y_t\|_2^{9/5}\|Y_t\|_{1,2}^{3/5})\\
& \leq C(\|\Delta Y_t\|_2^2 + \|Y_t\|_{1,2}^4\|X_t\|_{1,2}^2 + \|\Delta X_t\|_2^2 + \|X_t\|_{1,2}^6 + \|Y_t\|_{1,2}^6  + \|\Delta Y_t\|_2^2\\
& \leq C_3 \|\Delta
Y_t\|_2^2 + C_4\|\Delta X_t\|_2^{2} + C_5(\|Y_t\|_{1,2}^6 +
\|X_t\|_{1,2}^6).
\end{aligned}
\end{equation*} 
\begin{proposition}
Assume that $a \in L^2(\Omega, L^2(0,T; \cW^{2,2}(\mathbb{T}^2)))$. Then there exists a constant $C=C(R)$ such that 
\begin{equation*}
    \displaystyle\int_0^t f_R(a_s)^2\|\mathcal{L}_{u_s}v_s\|_2^2 ds \leq C(R) \left [ \displaystyle\int_0^t \|v_s\|_{2,2}^2 ds+t\right]
\end{equation*}
\begin{equation*}
   \displaystyle\int_0^t f_R(a_s)^2\|\nabla \cdot(h_sv_s)\|_2^2 ds \leq C(R) \left [ \displaystyle\int_0^t \|v_s\|_{2,2}^2 ds + \displaystyle\int_0^t \|h_s\|_{2,2}^2 ds + 1\right].
\end{equation*}
\end{proposition}
\noindent\textbf{\textit{Proof.}}
We have 
\begin{equation*}
    \|\mathcal{L}_{u_s}v_s\|_2^2 \leq \|\mathcal{L}_{v_s}v_s\|_2^2 + \left( \frac{R}{\epsilon}\right)^2\|\nabla v_s\|_2^2
\end{equation*}
\begin{equation*}
    f_R(a_s)^2\|\mathcal{L}_{u_s}v_s\|_2^2 \leq f_R(a_s)^2\|\mathcal{L}_{v_s}v_s\|_2^2 + \left( \frac{R}{\epsilon}\right)^2(R+1)^2
\end{equation*}
and
\begin{equation*}
    \begin{aligned}
    \|\mathcal{L}_{v_s}v_s\|_2^2 & \leq \|v_s\|_4^2\|\nabla v_s\|_4^2 
    \leq C\|v_s\|_2\|\nabla v_s\|_2 \|\nabla v_s\|_2\|\nabla\nabla v_s\|_2. 
    \end{aligned}
\end{equation*}
Hence
\begin{equation*}
    f_R(a_s)^2\|\mathcal{L}_{v_s}v_s\|_2^2 \leq C(R+1)^3\|v_s\|_{2,2}
\end{equation*}
and
\begin{equation*}
   \begin{aligned}
    \displaystyle\int_0^t f_R(a_s)^2\|\mathcal{L}_{u_s}v_s\|_2^2 ds &\leq C(R+1)^3\displaystyle\int_0^t\|v_s\|_{2,2}ds + \left( \frac{R}{\epsilon}\right)^2(R+1)^3t \\
    & \leq C(R)\displaystyle\int_0^t(\|v_s\|_{2,2}^2 + 1)ds.
   \end{aligned}
\end{equation*}
Similarly 
\begin{equation*}
\begin{aligned}
    \|\nabla \cdot (h_su_s))\|_2^2 & \leq \|\nabla h_s\|_4^2\|u_s\|_4^2 + \|h_s\|_4^2\|\nabla u_s\|_4^2 \\
    & \leq C(\|\nabla h_s\|_2\|\nabla \nabla h_s\|_2\|u_s\|_2\|\nabla v_s\|_2 + \|h_s\|_2\|\nabla h_s\|_2\|\nabla u_s\|_2\|\nabla \nabla u_s\|_2)
\end{aligned}
\end{equation*}
and
\begin{equation*}
    \begin{aligned}
    f_R(a_s)\|\nabla \cdot (h_su_s))\|_2^2 \leq C(R+1)^3(\|h_s\|_{2,2} + \|u_s\|_{2,2}),
    \end{aligned}
\end{equation*}
therefore 
\begin{equation*}
   \displaystyle\int_0^t f_R(a_s)^2\|\nabla \cdot(h_sv_s)\|_2^2 ds \leq C(R) \left [ \displaystyle\int_0^t \|v_s\|_{2,2}^2 ds + \displaystyle\int_0^t \|h_s\|_{2,2}^2 ds + 1\right].
\end{equation*}

\noindent\textbf{\textit{Proof of Lemma \ref{lemma:globsol1}}}
One can write 
\begin{equation*}
\begin{aligned}
d\|a_{t\wedge\tau_R}\|_{1,2}^2  + 2\gamma\|a
_{t\wedge\tau_R}\|_{2,2}^2dt  =& 2\left(\langle v_{t\wedge\tau_R} +  \Delta v_{t\wedge\tau_R}, Q_{v_{t\wedge\tau_R}}\rangle + \langle h_{t\wedge\tau_R} +  \Delta h_{t\wedge\tau_R}, Q_{h_{t\wedge\tau_R}}\rangle \right)dt\\
& - 2 \langle v_{t\wedge\tau_R} + \Delta v_{\wedge\tau_R} ,fk\times v_{t\wedge\tau_R} + g\nabla p_{t\wedge\tau_R}\rangle dt \\
 &+ \displaystyle\sum_{i=1}^{\infty} \left(\langle (\mathcal{L}_i+\mathcal{A}_i)v_{t\wedge\tau_R}, (\mathcal{L}_i+\mathcal{A}_i) v_{t\wedge\tau_R} \rangle  +  \langle \mathcal{L}_i h_{t\wedge\tau_R}, \mathcal{L}_i h_{t\wedge\tau_R}\rangle\right)dt \\
 & + \displaystyle\sum_{i=1}^{\infty} \left( \langle \Delta v_{t\wedge\tau_R},(\mathcal{L}_i+\mathcal{A}_i)^2 v_{t\wedge\tau_R}\rangle  +  \langle \Delta h_{t\wedge\tau_R}, \mathcal{L}_i^2 h_{t\wedge\tau_R}\rangle \right)dt \\
& + 2\displaystyle\sum_{i=1}^{\infty}\left( \langle v_{t\wedge\tau_R}, (\mathcal{L}_i+\mathcal{A}_i) v_{t\wedge\tau_R}\rangle +  \langle h_{t\wedge\tau_R}, \mathcal{L}_i h_{t\wedge\tau_R}\rangle \right)dW_t^i\\
& + 2\displaystyle\sum_{i=1}^{\infty}\left(\langle \Delta v_{t\wedge\tau_R}, (\mathcal{L}_i+\mathcal{A}_i) v_{t\wedge\tau_R}\rangle +  \langle \Delta h_{t\wedge\tau_R}, \mathcal{L}_i h_{t\wedge\tau_R} \rangle \right)
dW_t^i.  
\end{aligned}
\end{equation*}
Define
\begin{equation*}
\begin{aligned}
    \tilde{F}(a_s) :=& 2\left(\langle v_{s} +  \Delta v_{s}, Q_{v_{s}}\rangle + \langle h_{s} +  \Delta h_{s}, Q_{h_{s}}\rangle - \langle v_{s} + \Delta v_{s} ,fk\times v_{s} + g\nabla p_{s}\rangle \right) \\
 &+ \displaystyle\sum_{i=1}^{\infty} \left(\langle (\mathcal{L}_i+\mathcal{A}_i)v_{s}, (\mathcal{L}_i+\mathcal{A}_i) v_{s} \rangle  +  \langle \mathcal{L}_i h_{s}, \mathcal{L}_i h_{s}\rangle + \langle \Delta v_{s},(\mathcal{L}_i+\mathcal{A}_i)^2 v_{s}\rangle  +  \langle \Delta h_{s}, \mathcal{L}_i^2 h_{s}\rangle \right)
\end{aligned}
\end{equation*}
\begin{equation*}
\begin{aligned}
    \tilde{G}_i :=& 2\displaystyle\sum_{i=1}^{\infty}\left( \langle v_{s}, (\mathcal{L}_i+\mathcal{A}_i) v_{s}\rangle +  \langle h_{s}, \mathcal{L}_i h_{s}\rangle + \langle \Delta v_{s}, (\mathcal{L}_i+\mathcal{A}_i) v_{s}\rangle +  \langle \Delta h_{s}, \mathcal{L}_i h_{s} \rangle \right). 
\end{aligned} 
\end{equation*}
Then by Lemma \ref{abstractlemmaexistence} we have that 
\begin{equation}\label{contr1}
\|\tilde{F}(a_s)\|_{1,2}^2 \leq C\|a_s\|_{1,2}^6 - \zeta \|a_s\|_{2,2}^2 \leq C\|a_s\|_{1,2}^6 - \zeta \|a_s\|_{1,2}^2
\end{equation}
since we can choose $\zeta < 2\gamma$ such that \eqref{contr1} holds.  \\

\noindent The following result is introduced in Theorem 2, pp. 133, in \cite{Rozovskii}, for the $d$-dimensional domain $\mathbb{R}^d$. We rewrite it here for the two-dimensional torus $\mathbb{T}^2$: 
\begin{theorem}\label{theoremrozovskii}
Suppose that the following conditions hold true:
\begin{itemize}
\item[a.] $2\sigma^{ij}(x)\alpha^i\alpha^j - \displaystyle\sum_{i=1}^2|\sigma^{ij}(x)\alpha^i|^2 \geq b\displaystyle\sum_{i=1}^2|\sigma^i|^2, \ \ \  \forall \alpha \in\mathbb{T}^2$, where $b>0$ is independent of $t, \omega, x, \alpha$.
\item[b.] The functions $a^{ij}, b^i, c, \sigma^{il}, h^l$ with $i,j,l=1,2$ are differentiable in the spatial variable $x$ up to order $k$, for all $t, x, \omega$. Moreover, they are uniformly bounded (with respect to all variables) together with their derivatives, by a constant $C$. 
\item[c.] $u_0 \in L^2(\Omega, \mathcal{W}^{k,2}(\mathbb{T}^2))$, $f\in L^2([0,T]\times\Omega;\mathcal{W}^{k-1}(\mathbb{T}^2))$, $g^l \in L^2([0,T]\times\Omega;\mathcal{W}^{k,2}(\mathbb{T}^2)), l=1,2$.
\end{itemize}
Then the generalized solution $u$ of the problem 
\begin{equation*}
    \begin{aligned}
    du(t,x,\omega) & = \left[ (a^{ij}(t,x,\omega)u_i(t,x,\omega))_j + b^i(t,x,\omega)u_i(t,x,\omega) + c(t,x,\omega)u(t,x,\omega) + f(t,x,\omega)\right]dt\\
    &+ \left[\sigma^{ij}(t,x,\omega)u_i(t,x,\omega) + h^i(t,x,\omega)u(t,x,\omega) + g^l(t,x,\omega) \right]dW_t^l, \ \ \ (t,x,\omega) \in (0,T]\times\mathbb{T}^2\times\Omega, \\
    \end{aligned}
\end{equation*}
belongs to the class $L^2([0,T],\mathcal{W}^{k+1}(\mathbb{T}^2)) \cap C(0,T; \mathcal{W}^{k,2}(\mathbb{T}^2))$ and there exists $N = N(C,k,T) > 0$ such that 
\begin{equation*}
    \mathbb{E}\left[ \displaystyle\sup_{t\in[0,T]}\|u_t\|_{k,2}^2 + \displaystyle\int_0^T \|u_t\|_{k+1,2}^2dt\right] \leq N\mathbb{E}\left[ \|u_0\|_{k,2}^2 + \displaystyle\int_0^T\left(\|f_t\|_{k-1,2}^2 + \displaystyle\sum_{l=1}^2\|g_t^l\|_{k,2}^2 \right)dt\right].
\end{equation*}
\end{theorem}


\begin{thebibliography}{99}

\bibitem{Billingsley}
Billingsley, P., Convergence of Probability Measures, Second Edition, John Wiley \& Sons. (1999). 

\bibitem{BreschsDesjardinsMetivier}
Bresch D., Desjardins B., Métivier G. (2006) Recent Mathematical Results and Open Problems about Shallow Water Equations. In: Calgaro C., Coulombel JF., Goudon T. (eds) Analysis and Simulation of Fluid Dynamics, Advances in Mathematical Fluid Mechanics. Birkhäuser Basel, https://doi.org/10.1007/978-3-7643-7742-7\_2.

\bibitem{BreschsDesjardins1}
Bresch, D., Desjardins, B., Existence of Global Weak Solutions for a 2D Viscous Shallow Water Equations and Convergence to the Quasi-Geostrophic Model. Commun. Math. Phys. 238, 211–223 (2003). https://doi.org/10.1007/s00220-003-0859-8. 

\bibitem{BreschsDesjardins2}
Bresch, D., Desjardins, B., On the construction of approximate solutions for the 2D viscous shallow water model and for compressible Navier–Stokes models, Journal de Mathématiques Pures et Appliquées
Volume 86, Issue 4, October 2006, Pages 362-368, https://doi.org/10.1016/j.matpur.2006.06.005. 

\bibitem{Breschs2} 
Bresch, D., Desjardins, B., On the construction of approximate solutions for the 2D viscous shallow water model and for compressible NavierStokes models, J. Math. Pures Appl. 86 (2006) 362368.

\bibitem{Bui}
Bui, A. T., Existence and uniqueness of a classical solution of an initial boundary value
problem of the theory of shallow waters, SIAM J. Math. Anal. 12 (1981), 229-241.

\bibitem{ChenMiaoZhang}
Chen, Q., Miao, C., Zhang, Z., On the Well-Posedness for the Viscous Shallow Water Equations, SIAM J. Math. Anal., 40(2), 443–474, https://doi.org/10.1137/060660552. 


\bibitem{rotRSW1}
Cheng, B., Tadmor, E., Long-time Existence of Smooth Solutions for the Rapidly Rotating Shallow-Water and Euler Equations, SIAM J. Math. Anal. Vol. 39, No. 5, pp. 16681685 (2008).

\bibitem{CotterCrisanPan1}
Cotter, C. et. al., Numerically Modelling Stochastic Lie Transport in Fluid Dynamics (2018), available here: https://arxiv.org/abs/1801.09729.

\bibitem{CotterCrisanPan2}
Cotter, C. et. al., Modelling uncertainty using circulation-preserving stochas- tic transport noise in a 2-layer quasi-geostrophic model (2018) available here: https://arxiv.org/abs/1802.05711.

\bibitem{CL1}
Crisan, D., Lang, O., Well-posedness for a stochastic 2D Euler equation with transport noise, https://arxiv.org/abs/1907.00451 (2020).

\bibitem{CL2}
Crisan, D., Lang, O., Local well-posedness for the great lake equation with transport noise, Romanian Journal of Pure and Applied Mathematics, No 1 (2021).

\bibitem{Darryllecturenotes} 
Crisan, D., et. al., Mathematics Of Planet Earth: A Primer (Chapter 2), Advanced Textbooks In Mathematics (2017).

\bibitem{Temam2} 
Cyr, J., Nguyen, P., Temam, R., Stochastic one layer shallow water equations with Lévy noise, Discrete and Continuous Dynamical Systems, Series B, Volume 24, Number 8, August 2019.

\bibitem{DaPratoZabczyk} 
Da Prato, G., Zabczyk, J., Stochastic equations in infinite dimensions, Second Edition, Cambridge University Press, 2014, ISBN 978-1-107-05584-1 Hardback.

\bibitem{Kurtzprel} 
Ethier, S., Kurtz, T., Markov Processes - Characterization and Convergence, Wiley \& Sons, 1986, ISBN-I0 0-471-76986-X.


\bibitem{Holm2015} 
Holm, D., Variational principles for stochastic fluid dynamics, Proc. R.Soc.A 471:20140963 (2015).


\bibitem{HolmLuesink} 
Holm, D., Luesink, E., Stochastic wave-current interaction in stratified shallow water dy- namics, arXiv:1910.10627.

\bibitem{Kalnay} 
Kalnay, E., Atmoshperic Modeling, Data Assimilation and Predictability, Cambridge University Press (2003).

\bibitem{Kloeden}
Kloeden, P., Global existence of classical solutions in the dissipative shallow water
equations, SIAM J. Math. Anal. 16 (1985) , 301-315.

\bibitem{KurtzProtter}
Kurtz, T.G., Protter, P.E. Weak convergence of stochastic integrals and differential equations II: Infinite dimensional case In: Talay D., Tubaro L. (eds) Probabilistic Models for Nonlinear Partial Differential Equations. Lecture Notes in Mathematics, vol 1627. Springer, Berlin, Heidelberg.


\bibitem{OanaPhD} 
Lang, O., Nonlinear stochastic transport partial differential equations: well-posedness and data assimilation, Phd Thesis (2020), https://doi.org/10.25560/89816.   

\bibitem{LiHongZhu}
Li, J., Hong, P. \& Zhu, W. Ill-posedness for the 2D viscous shallow water equations in the critical Besov spaces. J. Evol. Equ. 20, 1287–1299 (2020). https://doi.org/10.1007/s00028-019-00556-y. 

\bibitem{Temam1} 
Link, J., Nguyen, P., Temam, R., Local martingale solutions to the stochastic one layer shallow water equations, J. Math. Anal. Appl. 448 (2017) 93139.


\bibitem{LiuYin2}
Liu, Y., Yin, Z., Global existence and local well-posedness of the 2D viscous shallow water system in Sobolev spaces, Applicable Analysis, Vol. 95, No 1 (2016), Pages 78-96, Taylor \& Francis, doi:10.1080/00036811.2014.998205. 

\bibitem{LiuYin1}
Liu, Y., Yin, Z., Global existence and well-posedness of the 2D viscous shallow water system in Sobolev spaces with low regularity, Journal of Mathematical Analysis and Applications
Volume 438, Issue 1, 1 June 2016, Pages 14-28, https://doi.org/10.1016/j.jmaa.2016.01.046. 

\bibitem{LiuYin0}
Liu, Y., Yin, Z., Global existence and well-posedness of the 2D viscous shallow water system in Besov spaces, Nonlinear Analysis: Real World Applications
Volume 24, August 2015, Pages 1-17, https://doi.org/10.1016/j.nonrwa.2014.12.005. 

\bibitem{rotRSW2}
Liu, H., Tadmor, E., Rotation prevents finite-time breakdown, Physica D: Nonlinear Phenomena Volume 188, Issues 3–4, 1 February 2004, Pages 262-276, https://doi.org/10.1016/j.physd.2003.07.006. 

\bibitem{Matsumura}
 Matsumura, A., Nishida, T., Initial boundary value problems for the equations of
motion of general fluids, in \textit{Computing Methods in Applied Sciences and Engineering}
(R. Glowinski and F. Lions, Eds.), Vol. 5, pp. 389-406, North-Holland, Amsterdam (1982). 

\bibitem{Orenga} 
Orenga, P., Un théorème d’existence de solutions d’un probl`eme de shallow water, Arch.
Rational Mech. Anal. 130 (1995) 183-204. 9 Springer-Verlag (1995).

\bibitem{Roeckner} 
Röckner, M., Schmuland, B., Zhang, X., Yamada-Watanabe theorem for stochastic evolution equations in infinite dimensions, Condensed Matter Physics 2008, Vol. 11, No 2(54), pp. 247259

\bibitem{Rozovskii}
Rozovskii, R. L., Stochastic Evolution Systems, Kluwer Academic Publishers, 1990, 978- 0792300373.

\bibitem{SIM} 
Simon, J., \textit{Compact sets in the space $L^{p}(0,T;B)$}, Ann. Mat. Pura Appl. 146 (1987) pp 65-96.

\bibitem{SC2020} 
Street, O., Crisan, D., Semi-martingale driven variational principles, arXiv:2001.10105 (2020).


\bibitem{Sundbye}
 Sundbye. L., Global existence for the Dirichlet problem for the viscous shallow water equations.
J. Math. Anal. Appl. 202(1), 236–258 (1996). 

\bibitem{Vallis} 
Vallis, G. K., Climate and the Oceans, Princeton Primers in Climate (2012).

\bibitem{LeeuwenCrisanPotthastLang}
van Leeuwen, P. J., Lang, O., Crisan, D., Potthast, R., Data assimilation for SALT - Rotating Shallow Water Models (in preparation). 

\bibitem{Zeitlinbook} 
Zeitlin, V., Geophysical fluid dynamics: understanding (almost) everything with rotating shallow water models, Oxford University Press (2018).
\end{thebibliography}
\end{document}